\documentclass[12pt,letterpaper]{amsart}

\usepackage{fullpage}

\usepackage[light,none,bottom,german]{draftcopy}
\newcommand{\submod}[2]{{#1   ({#2})}}

\newcount\timehh\newcount\timemm
\timehh=\time 
\divide\timehh by 60 \timemm=\time
\usepackage{amssymb,amscd}

\newcommand{\m}{{\mathfrak m}}
\newcommand{\p}{{\mathfrak p}}

\newcommand{\F}{{\mathbf F}}

\newcommand{\Q}{{\mathbf Q}}
\newcommand{\R}{{\mathbf R}}
\newcommand{\Z}{{\mathbf Z}}

\newcommand{\vsp}{\vspace{8pt}}

\newcommand{\vtsp}{\vspace{16pt}}
\newcommand{\rarr}{{\rightarrow}}

\newcommand{\mymod}[2]{{ #1 \: (\bmod \: {#2})}}

\newtheorem*{thm}{Main Theorem}
\newtheorem{lem}{Lemma}

\newtheorem{prop}{Proposition}
\newtheorem{cor}{Corollary}

\newtheorem*{conj}{Conjecture}

\theoremstyle{remark}
\newtheorem{remark}{Remark}
\newtheorem{quest}{Question}

\newtheorem*{ack}{Acknowledgment}

\newcommand{\ov}[1]{{\overline{{#1}}}}

\begin{document}

\title{Moments of the rank of elliptic curves}

\author{Siman Wong
  {
    \protect \protect\sc\today\ -- 
    \ifnum\timehh<10 0\fi\number\timehh\,:\,\ifnum\timemm<10 0\fi\number\timemm
    \protect \, \, \protect \bf DRAFT
  }
}

\address{Department of Mathematics \& Statistics, University of Massachusetts.
	Amherst, MA 01003-4515 USA}

\email{siman@math.umass.edu}



\date{\today}


\begin{abstract}
Fix an elliptic curve $E/\Q$, and assume the generalized Riemann hypothesis
for the $L$-function
$
L(E_D, s)
$
for every quadratic twist $E_D$ of $E$ by $D\in\Z$.  We combine Weil's
explicit formula with techniques of Heath-Brown to derive an asymptotic
upper bound for the weighted moments of the analytic rank of $E_D$.  It follows
from this that, for any unbounded increasing function $f$ on $\R$,
the analytic rank and (assuming in addition the Birch-Swinnerton-Dyer
conjecture)
the number of integral points of $E_D$ are less than $f(D)$
for almost all $D$.
We also derive an upper bound for the density of low-lying zeros of $L(E_D, s)$
which is compatible with the random matrix models of Katz and Sarnak.
%
%

\end{abstract}

\maketitle

\tableofcontents

\section{Introduction}
\newcounter{property}

Let $E$ be an elliptic curve over $\Q$.  The Birch-Swinnerton-Dyer conjecture
predicts that
$$
r_{\text{mw}}(E) := \text{ the rank of the Mordell-Weil group of $E/\Q$}
$$
is equal to the \textit{analytic rank}
$$
r_{\text{an}}(E) := \text{ the order at $s=1$ of the $L$-function $L(E, s)$.}
$$
This implies in particular the Parity Conjecture:
$$
w(E) = (-1)^{r_{\text{mw}}(E)},
$$
where $w(E)$ denotes the sign of the functional equation of $L(E, s)$.  Denote
by $N_E$ the conductor of $E/\Q$, and by $E_D$ the quadratic twist of $E$ by
an integer $D$.  If $D$ is square-free and is prime to $2N_E$, we have the
relation \cite{mtt}
$$
w(E_D) = w(E) \chi_D(-N_E),
$$
where $\chi_D$ denotes the quadratic character associated to $\Q(\sqrt{D})$.
Thus among the square-free integers $D$ prime to $2N_E$, the Parity Conjecture
implies that
half of the twists $E_D$ have odd Mordell-Weil rank, and the other half, even.
Early experimental investigations (\cite{brumer-mcguinness},
\cite{zagier-kramarz})
suggest that a positive portion of the quadratic (resp.~cubic) twists have
rank $\ge 2$. 
On the other hand, the random matrix models of Katz and Sarnak
(\cite[$\S 4$ and $5$]{katz-sarnak:survey}, \cite[p.~9-10]{katz}),
which presupposes the  Generalized Riemann hypothesis (GRH), predicts that
half of the twists should have analytic rank $0$, and the other half,
analytic rank $1$, whence the average analytic rank over all twists should be
$1/2$.
See \cite{rubin} for a recent survey on ranks of elliptic curves, and 
\cite{katz-sarnak:survey} on random matrix theory.

\vsp

Goldfeld seems to have been the first person to investigate the average rank of elliptic
curves in a quadratic twist family.  His main tool is Weil's explicit formula.
For the rest of this paper $F$ denotes the triangle function
$$
F(x) = \max(0, 1-|x|).
$$
The explicit formula says that 
the sum over powers of traces of Frobenius of $E_D$, weighted by $F$,
is essentially equal to a sum of the Mellin transform of $F$ extended
over the non-trivial zeros of $L(E_D, s)$.  Under GRH, each term of this latter
sum is non-negative.  Since
$
r_{\text{an}}(E_D)
$
is the order of $L(E_D, s)$ at $s=1$, to bound the average analytic
rank we are led to study the average of the non-Archimedean side of
the twisted explicit formula.  In this way, Goldfeld \cite{goldfeld}
shows that under GRH,  for $x\gg_{E, \epsilon} 1$ we have
\begin{equation}
\sum_{|D|<x } r_\text{an}(E_D) \le (3.25 + \epsilon) \sum_{|D|<x} 1.
	\label{gold}
\end{equation}
He also points out that any improvement of the constant $3.25$ to a number
strictly less than $2$ would imply that a positive portion of the twists
would have analytic rank $0$, a statement which at present has been
proved unconditionally only for special classes of $E$.
In his unpublished manuscript, Heath-Brown \cite{heath-brown} makes a
major breakthrough by improving Goldfeld's constant, also under GRH, from
$3.25$ to $1.5$, and
with $D$ restricted to twists with the same root number.
This implies that under GRH, a positive
portion of the twists of $E$ have rank $0$ and $1$, respectively.
This improvement is a result of better control over the non-Archimedean
side of the twisted explicit formula, so Heath-Brown's upper bounds are in fact
upper bounds for the average of the Archimedean side.  By keeping track of
the contribution from all the non-trivial zeros and not just $s=1$, we can
apply Heath-Brown's technique to get an \textit{asymptotic formula}
for all moments of the twisted explicit formula.

\vsp

For the rest of this paper, the constants involved in any $O$, $o$ and $\ll$
expressions are with respect
to the variable $x$ \textit{only} and  depend only
on those parameters printed as subscripts next to these symbols.  In
particular, any unadorned $O$, $o$ and $\ll$ constants are absolute.

\vsp

\begin{thm}
          \label{thm:asymptotic}
Fix a non-negative, thrice continuously differentiable function $W$
compactly supported on $(0, 1)$ or $(-1, 0)$.
Fix an elliptic curve $E/\Q$, and
assume the GRH for every $L(E_D, s)$.
For any positive integer
$
k = o_E(\log\log\log x)
$,
as $1+i\tau_D$ runs through the non-trivial zeros of $L(E_D, s)$ with
$
\tau_D\not=0
$
we have
\begin{eqnarray*}
\lefteqn{\sum_D
   \Bigl[
     r_{\text{\rm an}}(E_D)
     +
     \sum_{\tau_D\not=0}
       \Bigl(
         \frac{\sin(\tau_D (\log x) / 2)}{\tau_D (\log x) /2}
       \Bigr)^2
   \Bigr]^k
   \:
   W\Bigl(\frac{D}{x^{k/2}\log^{2k+2}x}\Bigr)}
\\
&\le&
\frac{1}{2}
\Bigl[
  \Bigl(
    k + \frac{1}{2} + \frac{1}{\sqrt{3}} 
  \Bigr)^k
+
  \Bigl(
    k + \frac{1}{2} - \frac{1}{\sqrt{3}} 
  \Bigr)^k
+
  o_{E, W}(1)
\Bigr]
\sum_D
   W\Bigl( \frac{D}{x^{k/2}\log^{2k+2}x} \Bigr).
\end{eqnarray*}
\end{thm}

\vsp

Note that the Main Theorem is effective with respect to $k$ so long as
$
k = o_E(\log\log\log x)
$.
This allows us to deduce the following result (cf.~$\S \ref{sec:cor}$).

\vsp

\begin{cor}
       \label{cor:slow}
Let $f$ be an unbounded increasing function on $\R$.  Fix an elliptic curve
$E/\Q$, and assume the GRH for every $L(E_D, s)$.  Then the set of integers
$D$ for which 
$
r_{\text{an}}(E_D) > f(D)
$
has density zero.
\end{cor}

\vsp

Conjectures of Lang (and others) giving height bounds for rational and
integral points on elliptic curves suggest that `most' elliptic curves have no
integral
  points.\footnote{I
  would like to thank Professor Silverman for bringing this to my
  attention.}
Thanks
to Corollary \ref{cor:slow} and the work of Silverman, we can make this
precise for quadratic twist families.
Let 
\begin{equation}
E: y^2 = x^3 + Ax + B        \label{model}
\end{equation}
be a quasi-minimal model for $E/\Q$ (i.e.~$|4A^3 + 27B^2|$ is minimal subject
to $A, B\in\Z$).   Silverman \cite[Theorem A]{joe} shows that there exists an
absolute constant
$
\kappa
$
such that, if the $j$-invariant of $E/\Q$ is non-integral for $\le \delta$
primes, then
\begin{equation}
\Bigl[
  \begin{array}{ll}
    \text{the number of $S$-integral points}
    \\
    \text{on the quasi-minimal model (\ref{model})}
  \end{array}
\Bigr]
\le
\kappa^{(1+r_{\text{mw}}(E))(1+\delta) + \#S}.       \label{joe}
\end{equation}
Since (\ref{model}) is quasi-minimal for $E$, up to a bounded
power of $2$ and $3$ the Weierstrass
equation
\begin{equation}
y^2 = x^3 + AD^2 x + BD^3          \label{ed}
\end{equation}
is quasi-minimal for $E_D$ if $D$ is square-free.
Since the $j$-invariant is constant in a quadratic twist family, Silverman's
theorem plus Corollary \ref{cor:slow}  immediately yields the following
\textit{conditional} result which makes precise for quadratic twist families the
heuristic above on integral points (note that $N_{E_D}\ll_E D^2$).

\vsp

\begin{cor}
    \label{cor:integral}
Fix an elliptic curve $E/\Q$, and assume the GRH and the Birth-Swinnerton-Dyer
conjecture for every $L(E_D, s)$.  Then for any unbounded increasing function 
$f$ on $\R$, the set of integers $D$ for which the Weierstrass equation
{\rm(\ref{ed})} has more than $f(N_{E_D})$ integral points has density zero.
\qed
\end{cor}

\vsp

\begin{quest}
Brumer \cite{brumer} shows that the average analytic rank of all elliptic
curves over $\Q$, as ordered by their height, is $\le 2.3$.  Is there a higher
moment analog of this result?

Lang \cite[p.~140]{lang} conjectures that the number of integral points on a
quasi-minimal model of any $E/\Q$ should be bounded solely in terms of 
$
r_{\text{mw}}(E_D)
$.
Silverman \cite[p.~251]{joe1} poses the more precise conjecture that 
(\ref{joe}) should hold for all $E$ with no $\delta$-dependence.  Silverman's
conjecture
plus a higher moment analog of Brumer's theorem should allow us to extend
Corollary \ref{cor:integral} uniformly to all elliptic curves over $\Q$.
\end{quest}

\vsp

The two Corollaries above exploit the effectiveness of the Main Theorem with
respect to $k$.  We now investigate consequences of the Main Theorem for
fixed $k$.
First, we fix a number $R>0$ and set $k = [R/e]+1$ to obtain the following
weighted upper bound on
the density of large rank twists.

\vsp

\begin{cor}
          \label{cor:rank}
Fix an elliptic curve $E/\Q$, and the GRH for every $L(E_D, s)$.  Then for
any fixed $R>0$ and $x\gg_R 1$, we have
$$
\sum_{r_{\text{\rm an}}(E_D)\ge R}
     W\Bigl(\frac{D}{x}\Bigr)
\le
\frac{1/2+o_{E, W, R}(1)}
     {1.44467^R}
\sum_D
     W\Bigl(\frac{D}{x}\Bigr).
         \qed
$$
\end{cor}

\vsp

\begin{remark}
For $k=1$, the Main Theorem is essentially due to Heath-Brown
\cite{heath-brown}.  More precisely, denote by
$
\Delta_E(+)
$
and
$
\Delta_E(-)
$
the set of square-free integers $D$ prime to $N_E$ for which $L(E_D, s)$ have
root numbers $+1$ and $-1$, respectively.  Then Heath-Brown shows that
\begin{equation}
\sum_{D\in\Delta_E(\pm)} r_{\text{an}}(E_D) W\Bigl(\frac{D}{x}\Bigr)
\le
\Bigl(\frac{3}{2} + o_E(1) \Bigr)
\sum_D W\Bigl( \frac{D}{x} \Bigr).       \label{actual}
\end{equation}
It then follows that
\begin{eqnarray}
\sum_{\substack{D\in\Delta_E(+)
      \\r_{\text{an}}(E_D)=0}}
                                   W\Bigl(\frac{D}{x}\Bigr) 
&\ge&
\Bigl(
  \frac{1}{2} + o_E(1)
\Bigr)
\sum_{D\in\Delta_E(+)} W\Bigl(\frac{D}{x}\Bigr),
\label{plus}
\\
\sum_{\substack{D\in\Delta_E(-)
      \\r_{\text{an}}(E_D)=1}}
                                   W\Bigl(\frac{D}{x}\Bigr) 
&\ge&
\Bigl(
  \frac{3}{4} + o_E(1)
\Bigr)
\sum_{D\in\Delta_E(-)} W\Bigl(\frac{D}{x}\Bigr).
\label{minus}
\end{eqnarray}
The general outline of the proof of the Main Theorem follows that of
Heath-Brown; in particular, we make crucial use of his smooth averaging;
cf.~$\S \ref{sec:devote}$.
Our main contribution is in the handling of certain truncated multivariable
sums (Proposition \ref{prop:fundamental}) and in the arithmetic applications
(Corollary \ref{cor:slow} to \ref{cor:zero}).
In particular, for $k>1$ the Main Theorem (and hence Corollary \ref{cor:rank})
can also be refined to sum over $D\in\Delta_E(\pm)$ only; we can even drop
the condition $(D, N_E)=1$, at the cost of introducing tedious
congruence argument on $D$ in the proof of the Main Theorem.
Such refinements, however, do not improve the lower bounds (\ref{plus}) and
(\ref{minus}), 
 so we will not pursue these issues here.
\end{remark}

\vsp

From the proof of the Main Theorem we see that
$
x^{k/2}\log^{2k+2} x
$
can be replaced by $x^{k/2+\epsilon}$ for any $\epsilon>0$, provided that we
stipulate the $o(1)$-term on the right side be dependent upon $\epsilon$.  We
can then rewrite the Main Theorem in a more suggestive form:
\begin{eqnarray}
\lefteqn{\sum_D
  \Bigl[
    r_{\text{an}}(E_D)
    +
    \sum_{\tau_D\not=0}
      \Bigl(
        \frac{\sin ( \frac{\tau_D \log T}{k+\epsilon} )}
	     {\frac{\tau_D \log T}{k+\epsilon}}
      \Bigr)^2
  \Bigr]^k
  W\Bigl(
    \frac{D}{T}
  \Bigr)}            \nonumber
\\
&=&
\frac{1}{2}
\Bigl[
  \Bigl(
    k + \frac{1}{2} + \frac{1}{\sqrt{3}} + \epsilon
  \Bigr)^k
+
  \Bigl(
    k + \frac{1}{2} - \frac{1}{\sqrt{3}} + \epsilon
  \Bigr)^k
+
  o_{E, W, \epsilon}(1)
\Bigr]
\sum_D
	W\Bigl( \frac{D}{T} \Bigr).
                     \label{epsilon}
\end{eqnarray}
The number of non-trivial zeros $\rho$ of $L(E_D, s)$ with
$
|\text{im}(\rho)| < Y
$
is
$
\frac{Y\log Y|D|}{2\pi} + O_E(Y + \log |D|)
$.
Thus
$$
\sum_{|\tau_D| \gg_{E, k, \epsilon} 1}
      \Bigl(
        \frac{\sin ( \frac{\tau_D \log |D|}{k+\epsilon} )}
	     {\frac{\tau_D \log |D|}{k+\epsilon}}
      \Bigr)^2
\sim
(k+\epsilon)
\sum_{|\tau_D| \gg_{E, k, \epsilon} 1}
      \Bigl(
        \frac{\sin ( \tau_D \log |D| )}
	     {\tau_D \log |D| )}
      \Bigr)^2.
$$
This \textit{suggests} that \textit{if} the low-lying zeros of $L(E_D, s)$ are
uniformly
distributed as $D$ varies, then removing the factor
$
k+\epsilon
$
from the $\tau_D$-sum in (\ref{epsilon}) \textit{should} result in scaling the right
side of (\ref{epsilon}) by a factor of
$
(k+\epsilon)^{-k}
$.
That would mean almost all twists of $E_D$ would have analytic rank
$
\le 1+\epsilon
$.

\begin{quest}
Can we make precise this heuristic argument?  Specifically, does random
matrix theory provide the proper framework within which to formulate
the type of uniform distribution statement required here?
\end{quest}

The factor $k+\epsilon$ in the $\tau_D$-sum is due to the
fact that the asymptotic formula in (\ref{epsilon}) sums over 
$
|D|  \ll_W x^{k/2+\epsilon}
$.
If we can prove a similar formula -- even just an upper bound -- by summing
over
$
|D|  \ll_W x^{\alpha}
$
for some fixed $\alpha$, uniformly for infinitely many $k$, then we would be
able to prove that almost all $E_D$ have analytic rank
$
\le 2\alpha + 1
$.
The reason we need to take such a long sum is to ensure that the main term
dominates the error term (\ref{become}).  Now, our argument leading up to
(\ref{become}) is essentially optimal, except in one step where we estimate
a \text{difference} of two terms by bounding each term; cf. Remark
\ref{remark:failure}.

\begin{quest}
Can we improve this error term (\ref{become})?
\end{quest}

\vsp

Corollary \ref{cor:rank}  gives an upper bound for the weighted average of the
multiplicity of the (potential) zero at $s=1$ of $L(E_D, s)$.  This
argument can be extended to count non-trivial zeros of bounded height.
We begin with some notation.  If $E_D$
is an even twist, then under GRH the non-trivial zeros of
$
L(E_D, s)
$
come in complex conjugate pairs $1+i\gamma_{E_D, j}$ with
$
0 \le \gamma_{E_D, 1} \le \gamma_{E_D, 2} \le \cdots
$.
If $E_D$ is an odd twist, then $L(E_D, s)$ has a zero at $s=1$; we label the
\textit{remaining} zeros as $1+i\gamma_{E_D, j}$ with
$
0 \le \gamma_{E_D, 1} \le \gamma_{E_D, 2} \le \cdots
$.
Finally, regardless of the parity of $E_D$, define
$$
\tilde{\gamma}_{E_D, j} = \gamma_{E_D, j} (\log N_{E_D}) / 2\pi.
$$
Since
$
\bigl( \sin(\frac{x}{2})/\frac{x}{2} \bigr)^2
$
is decreasing for $0<x< 2\pi$, if for some $|D|\gg_E 1$ we have
$
\tilde{\gamma}_{E_D, 3k} < 1/2\pi
$,
then for this $D$ and for every $j\le 3k$,
$$
\Bigl(
  \frac{\sin(\tau (\log |D|)/2)}{\tau(\log |D|)/2}
\Bigr)^2
>
\Bigl(
  \frac{\sin(1/2)}{1/2}
\Bigr)^2
=
0.9193953884.
$$
Invoke the Main Theorem and we get

\begin{cor}
         \label{cor:zero}
Fix an elliptic curve $E/\Q$, and assume the GRH for every $L(E_D, s)$.
For any integer $k>0$ and $x\gg_k 1$, we have
$$
\sum_{\tilde{\gamma}_{E_D, 3k} < 1/2\pi}
     W\Bigl(\frac{D}{x}\Bigr)
\le
\frac{1+o_{E, W, k}(1)}
     {1.402408^k}
\sum_D
     W\Bigl( \frac{D}{x} \Bigr).       \qed
$$
\end{cor}

\vsp

To put this result into context, recall that random matrix theory
\cite[$\S 6.9, \S 7.5.5$]{katz-sarnak} furnishes a family of probability measures
$
v(+, j), v(-, j)
$
on $\R$, $j=1, 2, \ldots$, with respect to which Katz and Sarnak formulate
the following conjecture.

\vsp

\begin{conj}[Katz-Sarnak]
For any integer $j\ge 1$ and any compactly supported complex-value function
$h$ on $\R$, 
$$
\sideset{}{'}\sum_{w(E_D)=+1}  h(\tilde{\gamma}_{E_D, j})
=
\Bigl(
  \sideset{}{'}\sum_{w(E_D)=+1} 1 + o_{E, h}(1)
\Bigr)
\int_\R h \cdot dv(+, j),
$$
where
$
\sideset{}{'}\sum_D
$
signifies that $D$ runs through all square-free integers $D$.
Similarly for $v(-, j)$.
\end{conj}

\vsp

As is pointed out in (\cite[p.~21]{katz-sarnak:survey}, \cite[p.~10]{katz}),
this Conjecture implies that almost all even (resp.~odd) twists of $E$ have
analytic rank $0$ (resp.~$1$).  By choosing $h$ to be supported on an
arbitrarily small
neighborhood of $0\in\R$, this Conjecture implies that for any fixed $j$ and
 any $\epsilon>0$, there exists 
$
\delta_j(\epsilon)>0
$
so that
\begin{itemize}
\item
$
\delta_j(\epsilon) \rarr 0
$
as $\epsilon\rarr 0$; and
\item
the set of square-free $D$ for which
$
\tilde{\gamma}_{E_D, j} < \epsilon
$
and $w(E_D)=1$, has density $< \delta_j(\epsilon)$.
\end{itemize}
In particular, for any $\epsilon>0$ the $\delta_j(\epsilon)$ (if they exist)
form a
non-increasing sequence that converges to $0$.
With respect to this formalism, Corollary \ref{cor:zero} can be viewed as
proving the existence of
$
\delta_j(1/2j)
$
under GRH (instead of the full random matrix theory), such that
$
\text{$\delta_{j}(1/2\pi) \rarr 0$ as $j\rarr\infty$}
$.
However, our present argument does not allow us to replace 
$
1.402408
$
with an arbitrarily large constant by replacing $1/2\pi$ with an arbitrarily
small number.

\vsp

\begin{remark}
The Main Theorem, and hence the Corollaries, readily extends to cubic and
higher order twists; cf.~Remark \ref{remark:cubic}.  We can also replace
$E$ by a newform.
\end{remark}

\vsp

\begin{ack}
I am indebted to Professor Heath-Brown for sending me a copy of his preprint
\cite{heath-brown}.  I would like to thank Professors Hajir, Hoffstein, Rosen
and Silverman for many useful discussions.
\end{ack}

\vsp

\section{Explicit formula}
    \label{sec:quad}

Fix a modular elliptic curve $E/\Q$ of conductor $N_E$.  Denote by $a_n(E)$
the $n$-th coefficient of $L(E, s)$.  For any prime
$
p\nmid N_E
$,
denote by 
$
\alpha_p(E)
$
and
$
\ov{\alpha}_p(E)
$
the eigenvalues of the Frobenius of $E/\F_p$.  Define
$$
c_n(E)
=
\left\{
  \begin{array}{lllll}
    \alpha_p(E)^m + \ov{\alpha}_p(E)^m
    &
    \text{if $n=p^m>1$ and $p\nmid N_E$;}
    \\
    a_p(E)^m
    &
    \text{if $n=p^m>1$ and $p| N_E$;}
    \\
    0
    &
    \text{otherwise.}
  \end{array}
\right.
$$
Note that $c_p(E) = a_p(E)$. 
For any $\lambda>0$, define
$
F_\lambda(x) = F(x/\lambda)
$.
Denote by $\Phi_\lambda(x)$ the Mellin transform of $F_\lambda$:
$$
\Phi_\lambda(u)
=
\int_{-\infty}^\infty  F_\lambda(x) e^{(u-1)x} dx.
$$
Note that if $s=1+it$ with $t\in\R$, then 
\begin{equation}
\Phi_\lambda(s)
=
\lambda
\Bigl(
   \frac{\sin(\lambda t/2)}{\lambda t/2}
\Bigr)^2.
      \label{phi}
\end{equation}
 As $\rho$ runs through the zeros 
$
\rho = \beta + i\gamma
$
of $L(E, s)$ with
$
0 < \beta < 2
$,
counted with multiplicity, Weil's explicit formula \cite[$\S$II.2]{mestre}
says that
\begin{eqnarray}
\lefteqn{\sum_\rho \Phi_\lambda(\rho)
:=
\lim_{z\rarr\infty} \sum_{|\rho|< z} \Phi_\lambda(\rho)}
        \nonumber
\\
&=&
  \log N_{E}
  -
  2\sum_{p^m>1} \frac{c_{p^m}(E)\log p}{p^m}
  F\Bigl(
      \frac{\log p^m}{\lambda}
   \Bigr)
  -
  2\log 2\pi
  -
  2\int_0^\infty
    \Bigl(
      \frac{F(t/\lambda)}{e^t -1} - \frac{1}{t e^t}
    \Bigr)
    dt.
    \label{weil}
\end{eqnarray}
Note that 
$
|c_{p^m}(E)| \le 2p^{m/2}
$.
Since $||F|| \le 1$, that means
$$
\sum_{\substack{p, m\\m\ge 3}}
      \frac{c_{p^m}(E)\log p}{p^m}
      F\Bigl( \frac{\log p^m}{\lambda}\Bigr)
\ll
\sum_{\substack{p, m\\m\ge 3}}
      \frac{\log p}{p^{m/2}}
\ll
\sum_{n>1}  \frac{\log n}{n^{3/2}}
\ll 
1.
$$
For $\lambda\ge 1$, the integral in (\ref{weil}) is bounded from above and
below by absolute constants, so the explicit formula now takes the form
$$
\sum_\rho \Phi_\lambda(\rho)
=
  \log N_{E}
  -
  2\sum_p
          \frac{c_{p}(E)\log p}{p}F\Bigl(\frac{\log p}{\lambda}\Bigr)
  -
  2\sum_p
          \frac{c_{p^2}(E)\log p}{p}F\Bigl(\frac{\log p^2}{\lambda}\Bigr)
  +
  O(1).
$$

Next, we study how the explicit formula behaves under quadratic twists.  If
$
p\nmid 2N_E D
$
(note that $2N_E$ and $D$ need not be coprime and $D$ need not be square-free),
then
$$
c_p(E_D) = a_p(E) \Bigl(\frac{D}{p}\Bigr),
\:\:
c_{p^2}(E_D) = c_{p^2}(E).
$$
Since $|| F || \le 1$,
\begin{eqnarray*}
\sum_{p | 2N_E D}
	F\Bigl(\frac{\log p}{\lambda}\Bigr)
	\frac{\log p}{p}
	\Bigl(
		c_p(E_D) - a_p(E)\Bigl(\frac{D}{p}\Bigr)
	\Bigr)
&\ll&
\sum_{p | 2N_E D} \frac{\log p}{\sqrt{p}},
\\
\sum_{p | 2N_E D}
	F\Bigl(\frac{\log p}{\lambda}\Bigr)
	\frac{\log p}{p}
	\Bigl(
		c_{p^2}(E_D) - c_{p^2}(E)
	\Bigr)
&\ll&
\sum_{p | 2N_E D} \frac{\log p}{p}.
\end{eqnarray*}
Since $\log p\ll p^{1/4}$, for $|D|\ge 2$ the right side of both expressions
above are
$$
\ll
\sum_{p|2 N_E D} p^{-1/4}
\ll
\sum_{p<\log(2N_E |D|)} p^{-1/4}
\ll_E
\log^{3/4} |D|.
$$
As $\rho_D$ runs through the zeros of $L(E_D, s)$ with $0<\text{Re}(\rho_D)<2$,
we now have
\begin{eqnarray*}
\sum_{\rho_D} \Phi_\lambda(\rho_D)
&=&
  \log N_{E_D}
  -
  2\sum_p
          \frac{c_p(E)\log p}{p}
	  \Bigl(\frac{D}{p}\Bigr)
	  F\Bigl(\frac{\log p}{\lambda}\Bigr)
\\
&&
\hspace{42pt}
  -
  \:
  2\sum_p
          \frac{c_{p^2}(E)\log p}{p^2}
	  F\Bigl(\frac{2\log p}{\lambda}\Bigr)
  +
  O(\log^{3/4} |D|).
\end{eqnarray*}

\begin{lem}
      \label{lem:rankin}
We have the estimates
\begin{eqnarray*}
\sum_p
         \frac{c_{p^2}(E)\log p}{p^2} F\Bigl( \frac{2\log p}{\lambda}\Bigr)
&=&
-\lambda/4 + o_E(\lambda),
\\
\sum_p
         \frac{a_p(E)^2\log^2 p}{p^2} F\Bigl( \frac{\log p}{\lambda}\Bigr)^2
&=&
\lambda^2 / 12 + o_E(\lambda^2).
\end{eqnarray*}
\end{lem}

\begin{proof}
If $p\nmid N_E$, then $c_{p^2}(E) = a_p(E)^2 - 2p$, so
$$
\sum_{p\nmid 2N_E} \frac{c_{p^2}(E)\log p}{p^s}
=
\sum_{p\nmid 2N_E} \frac{a_p(E)^2\log p}{p^s}
-
2\sum_{p\nmid 2N_E} \frac{\log p}{p^{s-1}}.
$$
Up to the bad primes and a term holomorphic for $Re(s)>3/2$, the two sums on
the right are $(-1)$ times the logarithmic derivative of, respectively,
the Rankin-Selberg $L$-function of the cusp form associated to $E$ with itself,
and $\zeta(s-1)$.  Each of the convolution $L$-function and $\zeta(s-1)$ has a
simple
pole at $s=2$.  Tauberian theorem then gives
$$
-\sum_{p<x} \frac{a_p(E)^2\log p}{p}
=
\sum_{p<x} \frac{c_{p^2}(E)\log p}{p}
=
-x + o_E(1).
$$
The Lemma now follows from partial summation.
\end{proof}

\vsp

Set $\lambda=\log x$ and define
$$
\beta_p
=
\frac{a_p(E)\log p}{p} F\Bigl( \frac{\log p}{\log x} \Bigr),
\:\:\:
X_k = x^{k/2}\log^{2k+2}x.
$$
In what follow, we will take $D$ so that $|D|\le X_k$.  From now on,
assume\footnote{We
  choose this $o$-bound for $k$ to simplify the exposition.  The optimal
  choice would be that which renders the $O$-term in Proposition
  \ref{prop:devoted} to be $o_E(\log x)$, but such refinements have no
  material impact on the arithmetic applications of the Main
  Theorem.}
\begin{equation}
k = o_E(\log\log\log x),               \label{k}
\end{equation}
whence
$
O_E(\log^{3/4}|D|) = o_E(\log x)
$.
Combine all these and recall
that $N_{E_D} \ll N_E D^2$, we now arrive at the final form of the explicit
formula for $E_D$:
\begin{eqnarray}
\sum_{\rho_D} \Phi_{\log x}(\rho_D)
&\le&
  \log (D^2)
  +
  (\log x)/2
  -
  2\sum_p \beta_p \Bigl( \frac{D}{p} \Bigr)
  +
  o_E(\log x).        \label{twisted}
\end{eqnarray}
We emphasize again that $D$ need not be coprime to $2N_E$ or square-free.

\vsp

\section{Moments of analytic rank}
     \label{sec:moments}

Define
\begin{eqnarray*}
f(x,D)
=
2\log |D| + (\log x)/2,
\:\:\:
R(x,D)
=
2\sum_p  
        \beta_p	\Bigl(\frac{D}{p}\Bigr).
\end{eqnarray*}
Let $W$ be a thrice continuously differentiable function with compact
support on $(1/2,1)$ or $(-1, -1/2)$.
The $k$-th moment of the twisted explicit formula, weighted by $W$, now becomes
\begin{eqnarray*}
\lefteqn{\sum_D
  \Bigl(\sum_{\rho_D}
         \Phi_{\log x}(\rho_D)\Bigr)^k
	 W\Bigl( \frac{D}{X_k}\Bigr)
\:\:
\le
\:\:
\sum_D
      \bigl( 2\log |D| + (\log x)/2 + o_E(\log x) \bigr)^k
      W\Bigl( \frac{D}{X_k} \Bigr)}
\\
&&
+
\sum_{r=1}^k
\Bigl(
  \!
  \begin{array}{l}k\\r\end{array}
  \!
\Bigr)
(-1)^r
\sum_D
f(x,D)^{k-r}  R(x,D)^r  W\Bigl(\frac{D}{X_k}\Bigr)
\\
&&
+
\sum_{r=1}^k
\Bigl(
  \!
  \begin{array}{l}k\\r\end{array}
  \!
\Bigr)
o_{E,k}\Bigl(
  \sum_{i=1}^{k-r}
  \Bigl(
    \!
    \begin{array}{c}k-r-i\\r\end{array}
    \!
  \Bigr)
  \log^i x
\sum_D
  f(x,D)^{k-r-i}  R(x,D)^r  W\Bigl(\frac{D}{X_k}\Bigr)
\Bigr)
\end{eqnarray*}
We begin by tackling the first of the three sums on the right.

\vsp

\begin{lem}
  \label{lem:firstterm}
For $l\ge 0$, we have
$$
\sum_D  f(x,D)^l W\Bigl(\frac{D}{X_k}\Bigr)
=
\Bigl( (k + 1/2)\log x + o_{E, W}(\log x)\Bigr)^l
\Bigl[
  \sum_D W\Bigl(\frac{D}{X_k}\Bigr) + o(X_k)
\Bigr].
$$
\end{lem}

\begin{proof}
Since $W(x)=0$ if $|x|\ge 1$, the sum in the Lemma extends over $|D|\le X_k$
only.  Thus with
$
X' := x^{k/2}
$,
from (\ref{k}) we see that
\begin{eqnarray*}
\lefteqn{ (k+1/2)\log x + o_E(\log x)\Bigr)^l
\sideset{}{'}\sum_{|D|> X'}
      W\Bigl(\frac{D}{X_k}\Bigr)
\ge
\sideset{}{'}\sum_{|D|> X'}
      f(x, D)^l      W\Bigl(\frac{D}{X_k}\Bigr)}
\\
&\hspace{90pt}
\ge&
\Bigl(  (k+1/2)\log x + o_E(\log x)\Bigr)^l
\sideset{}{'}\sum_{|D|> X'}
      W\Bigl(\frac{D}{X_k}\Bigr).
\end{eqnarray*}
The condition $|D| > X'$ can be dropped at the cost
of introducing a term
$$
\ll 
\Bigl( (k+1/2)\log x + o_E(\log x) \Bigr)^l
\sum_{|D|\le X'}
      W\Bigl(\frac{D}{X_k}\Bigr)
\ll_W
\Bigl( (k+1/2)\log x + o_E(\log x)\Bigr)^l
      x^{k/2},
$$
and the Lemma follows.
\end{proof}

\vsp

The rest of the paper is devoted to prove the following result.  The proof
of the Main Theorem makes use of the 
conditional estimate only; we state the unconditional result for
comparison.

\begin{prop}
                \label{prop:devoted}
For $r>0$, we have the estimate
\begin{eqnarray*}
\lefteqn{\sideset{}{'}\sum_D f(x,D)^i  R(x,D)^r  W\Bigl(\frac{D}{X_k}\Bigr)}
\\
&&
\left\{
  \renewcommand{\arraystretch}{1.6}
  \begin{array}{lll}
  =
  &
\displaystyle
  \Bigl(
      2\log X_k + \frac{\log x}{2} + o_{E, W}(\log x)
  \Bigr)^i
  (1+o(E))^{r/2}\log^r x 
  \sum_D
     W\Bigl( \frac{D}{X_k} \Bigr)
  \\
  &
  +
  \:
  O_{E, W}( 4^r r^3 x^{3r}(\log X_k+\log x)^{r+i}/ T^2)
          \hspace{100pt}
  \text{if $r$ is even,}
\\
  =
  &
   O_{E, W}( 4^r r^3 x^{3r}(\log X_k+\log x)^{r+i}/ T^2)
          \hspace{112pt}
  \text{if $r$ is odd.}
  \end{array}
\right.
  \renewcommand{\arraystretch}{1}
\end{eqnarray*}
If we assume the GRH for every $L(E_D, s)$, then the $O$-term can be improved
to
$$
O_{E,  W}( c_E^r r^{r+3} x^{r/2}(\log X_k+ \log x)^{r+i})
$$
for some constant $c_E$ depending on $E$ only.
\end{prop}

\vsp

Assuming the GRH-estimate, we then see that
\begin{eqnarray*}
\lefteqn{\frac{1}{\log^k x}
     \sum_D
     \Bigl(
       \sum_{\rho_D} \Phi_{\log x}(\rho_D)
     \Bigr)^k
     W\Bigl(\frac{D}{X_k}\Bigr)
\:\:
\le
\:\:
( k + 1/2 + o_E(1) )^k 
  \sum_D
     W\Bigl( \frac{D}{X_k} \Bigr)}
\\
&&
+
\sum_{\substack{r=1\\ \text{$r$ even}}}^k
\Bigl(
  \!
  \begin{array}{l}k\\r\end{array}
  \!
\Bigr)
(1+ o_E(1))^{r/2} 
(k + 1/2 + o_{E,k}(1) )^{k-r} 
(1/\sqrt{3})^r
  \sum_D
     W\Bigl( \frac{D}{X_k} \Bigr)
\\
&&
+
O_{E, W}( k^3 c_E^k x^{k/2}(\log X_k+ \log x)^{2k}).
\end{eqnarray*}
Recall (\ref{k}) and we see that this $O$-term is
$
o_{E, W}(X_k)
$.
To write
$
\sum_{\text{$r$ even}}
$
is to write
$
\frac{1}{2}\sum_{\text{all $r$}} ( 1 + (-1)^r)
$.
Expand the rest of the second line above accordingly and recall (\ref{phi}),
 we get
\begin{eqnarray}
\lefteqn{\sum_D
   \Bigl[
     r_{\text{\rm an}}(E_D)
     +
     \sum_{\tau_D\not=0}
       \Bigl(
         \frac{\sin(\tau_D (\log x) / 2)}{\tau_D (\log x) /2}
       \Bigr)^2
   \Bigr]^k
   \:
   W\Bigl(\frac{D}{X_k}\Bigr)}          \nonumber
\\
&\le&
\frac{1}{2}
\Bigl[
  \Bigl(
    k + \frac{1}{2} + \frac{1}{\sqrt{3}} 
  \Bigr)^k
+
  \Bigl(
    k + \frac{1}{2} - \frac{1}{\sqrt{3}} 
  \Bigr)^k
+
  o_{E, W}(1)
\Bigr]
\sum_D
   W\Bigl( \frac{D}{X_k} \Bigr),
\end{eqnarray}
and the Main Theorem follows.

\if 3\
{
To complete the proof of the Main Theorem it remains to sieve out the
contribution from the non-square-free $D$.

For any integer $n>0$,
\begin{eqnarray*}
&&
\sum_{D\equiv 0 \: (n^2)}
   \Bigl[
     r_{\text{\rm an}}(E_D)
     +
     \sum_{\tau_D\not=0}
       \Bigl(
         \frac{\sin(\tau_D (\log x) / 2)}{\tau_D (\log x) /2}
       \Bigr)^2
   \Bigr]^k
   \:
   W\Bigl(\frac{D}{X_k}\Bigr)
\\
&=&
\:\:\:
\sum_d
\:\:\:\:
   \Bigl[
     r_{\text{\rm an}}(E_d)
     +
     \sum_{\tau_d\not=0}
       \Bigl(
         \frac{\sin(\tau_d (\log x) / 2)}{\tau_d (\log x) /2}
       \Bigr)^2
   \Bigr]^k
   \:
   W\Bigl(\frac{d}{X_k/n^2}\Bigr)
\\
&=&
\frac{1}{2}
\Bigl[
  \Bigl(
    k + \frac{1}{2} + \frac{1}{\sqrt{3}} 
  \Bigr)^k
+
  \Bigl(
    k + \frac{1}{2} - \frac{1}{\sqrt{3}} 
  \Bigr)^k
+
  o_{E, k, W}(1)
\Bigr]
\sum_d
   W\Bigl( \frac{d}{X_k/n^2} \Bigr),
            \hspace{20pt}
	    \text{by (\ref{sqfree}).}
\end{eqnarray*}

\begin{lem}
            \label{lem:trap}
Denote by $\hat{W}$ the Fourier transform of $W$.  Then
$$
\sum_D W\Bigl(\frac{D}{x}\Bigr) = x\hat{W}(0) + O_W(1/x).
$$
\end{lem}

\begin{proof}
Set $x'=([x]+1)/x$.  The trapezoidal rule in numerical integration says that
$$
\Bigl|
  \sum_{|D|\le [x]+1} W\Bigl(\frac{D}{x}\Bigr)\frac{1}{x}
  -
  \int_{-x'}^{x'} W(t) dt
\Bigr|
\ll_W
x' \frac{1}{x^2}
\ll_W
1/x^2.
$$
Since $W$ has compact support on $(-1, 1)$, we are done.
\end{proof}

\vsp

Denote by $N_k(x)$ the set of square-free integers $0<n<X_k$ divisible only by
primes
$
\le \log x
$.
Then
\begin{eqnarray*}
&&
\sideset{}{'}\sum_D
   \Bigl[
     r_{\text{\rm an}}(E_D)
     +
     \sum_{\tau_D\not=0}
       \Bigl(
         \frac{\sin(\tau_D (\log x) / 2)}{\tau_D (\log x) /2}
       \Bigr)^2
   \Bigr]^k
   \:
   W\Bigl(\frac{D}{X_k}\Bigr)
\\
&=&
\sum_{n\in N_k(x)}
\mu(n)
\sum_d
   \Bigl[
     r_{\text{\rm an}}(E_d)
     +
     \sum_{\tau_d\not=0}
       \Bigl(
         \frac{\sin(\tau_d (\log x) / 2)}{\tau_d (\log x) /2}
       \Bigr)^2
   \Bigr]^k
   \:
   W\Bigl(\frac{d}{X_k/n^2}\Bigr)
\\
&&
+
O\Bigl(
\sum_{p>\log x}
\sum_d
   \Bigl[
     r_{\text{\rm an}}(E_d)
     +
     \sum_{\tau_d\not=0}
       \Bigl(
         \frac{\sin(\tau_d (\log x) / 2)}{\tau_d (\log x) /2}
       \Bigr)^2
   \Bigr]^k
   \:
   W\Bigl(\frac{d}{X_k/p^2}\Bigr)
\Bigr).
\end{eqnarray*}
Since
$
\sum_D W(D/T) \ll_W T
$,
the $O$-term is
\begin{eqnarray*}
&
\ll_{E, W}
&
\sum_{p> \log x} k^k X_k / p^2       \hspace{20pt}\text{by (\ref{sqfree})}
\\
&
\ll_{E, W}
&
k^k X_k /\log x
\\
&=&
o_{E, W}(X_k),       \hspace{40pt}\text{by (\ref{k})}
\end{eqnarray*}
while the main terms is
\begin{eqnarray*}
&&
\frac{1}{2}
\Bigl[
  \Bigl(
    k + \frac{1}{2} + \frac{1}{\sqrt{3}} 
  \Bigr)^k
+
  \Bigl(
    k + \frac{1}{2} - \frac{1}{\sqrt{3}} 
  \Bigr)^k
+
  o_{E, W}(1)
\Bigr]
\sum_{n\in N_k(x)}
\mu(n)
\sum_d W\Bigl(\frac{d}{X_k/n^2}\Bigr)
\\
&=&
\Bigl[
  \Bigl(
    k + \frac{1}{2} + \frac{1}{\sqrt{3}} 
  \Bigr)^k
+
  \Bigl(
    k + \frac{1}{2} - \frac{1}{\sqrt{3}} 
  \Bigr)^k
+
  o_{E, W}(1)
\Bigr]
\frac{X_k\hat{W}(0)}{2}
\Bigl[
  \sum_{n\in N_k(x)}  \frac{\mu(n)}{n^2}
  +
  o_{E, W}(1)
\Bigr]
\\
&=&
\Bigl[
  \Bigl(
    k + \frac{1}{2} + \frac{1}{\sqrt{3}} 
  \Bigr)^k
+
  \Bigl(
    k + \frac{1}{2} - \frac{1}{\sqrt{3}} 
  \Bigr)^k
+
  o_{E, W}(1)
\Bigr]
\frac{X_k\hat{W}(0)}{2}
\Bigl[
  \prod_{p<\log x}  \Bigl( 1 - \frac{1}{p^2} \Bigr)
  +
  o_{E, W}(1)
\Bigr]
\\
&=&
\Bigl[
  \Bigl(
    k + \frac{1}{2} + \frac{1}{\sqrt{3}} 
  \Bigr)^k
+
  \Bigl(
    k + \frac{1}{2} - \frac{1}{\sqrt{3}} 
  \Bigr)^k
+
  o_{E,  W}(1)
\Bigr]
\frac{X_k\hat{W}(0)}{2}
\Bigl[
   \frac{1}{\zeta(2)} + o_{E, W}(1)
\Bigr].
\end{eqnarray*}
Repeat this argument in reverse, we get
$$
\frac{X_k \hat{W}(0)}{2} \frac{1}{\zeta(2)}
=
\sideset{}{'}\sum_D
   W\Bigl(\frac{D}{X_k}\Bigr)
+
   o_{E, W}(X_k),
$$
and the Main Theorem for square-free $D$ follows.
\end{proof}
}
\fi

\vsp

\section{Proof of Corollary \ref{cor:slow}}
      \label{sec:cor}

Given any subset $S\subset\Z$, define its lower density to be the lim sup 
over all numbers $\sigma\ge 0$ such that
$$
\#\{ s\in S: |s|<x \} > \sigma x
  \hspace{10pt}
  \text{for all $x\gg_{S, \sigma} 1$}
$$
In particular, $S$ has density zero if and only if it has lower density zero.

\vsp

\begin{lem}
  \label{lem:lower}
With $W$ as in the Lemma, there exists a constant $\lambda_W>0$ depending on
$W$ only, such that for any subset $S\subset\Z$ with lower density
$
\sigma_S
$,
we have
$$
\sum_{s\in S} W\Bigl(\frac{s}{x}\Bigr) > \lambda_W \cdot \sigma_S x.
$$
\end{lem}

\begin{proof}
Without loss of generality, assume that $W$ is supported on $(0, 1)$.  Since
$W$ is continuous and since $W(0) = W(1)$, there exists an integer
$
n > 4/\sigma_S
$
such that for some 
$
0 < m < n-1
$
and some $w_0>0$, we have
$$
\text{$W(r)\ge w_0$
     for 
     $r\in N_n(m) := \{ t\in\R: m/n \le t\le (m+1)/n \}$.}
$$
For $0\le i< n$, set
$$
S_n(i, x) := \{ s\in S: i/n \le s/x \le (i+1)/n \}.
$$
By the definition of lower density, for
$
x\gg_S 1
$
we have
$
\sum_i \#S_n(i, x) \ge \sigma_S x
$.
That means 
$
\#S_n(j, x) \ge \sigma_S x/2n
$
for some $j\not=0, n-1$: otherwise
\begin{eqnarray*}
\sigma_S x
\:\:
\le
\:\:
\sum_i \#S_n(i, x)
&\le&
\sum_{i\not=0, n-1} \#S_n(i, x)
+
\#S_n(0, x)
+
\#S_n(n-1, x)
\\
&\le&
(n-2)\sigma_S x/2n + 2(x/n+ O(1))
\\
&=&
x(\sigma_S /2 - \sigma_S/n + 2/n) + O(1)
\\
&<&
(1-1/n)\sigma_S x + O(1),
\end{eqnarray*}
a contradiction.

Suppose this $j\ge m$; then for $s\in S_n(j, x)$,
\begin{eqnarray*}
\frac{m}{n}
\:\:
=
\:\:
\frac{j}{n} \frac{m}{j}
&\le&
\frac{s}{x}\frac{m}{j}
\\
&\le&
\frac{j+1}{n} \frac{m}{j}
\:\:
\le
\:\:
\frac{m+1}{n},
\end{eqnarray*}
whence
$$
\sum_{s\in S} W\Bigl(\frac{s}{xj/m}\Bigr)
\ge
\sum_{s\in S_n(j, x)} W\Bigl( \frac{s}{xj/m} \Bigr)
\ge
w_0 \frac{\sigma_S x}{n}
=
w_0 \frac{\sigma_S}{n}
\Bigl(
  \frac{xj}{m}
\Bigr)
\frac{m}{j}
\ge
w_0 \frac{\sigma_S}{n^2} 
\Bigl(
  \frac{xj}{m}
\Bigr).
$$
Next, suppose $j<m$.  Then for some $0\le l<m$ we have
$
\#S_l \ge \sigma_S x/ 2mn
$,
where
$$
S_l := \Bigl\{ s\in S: \frac{1}{n}\Bigl( j + \frac{l}{m}   \Bigr)
                       \le
		       \frac{s}{n}
                       \le
                       \frac{1}{n}\Bigr( j + \frac{l+1}{m} \Bigr) \Bigr\}.
$$
Multiplication by
$
\mu := (m+1)/(j + \frac{l+1}{m})
$
takes the interval 
$
[ j + \frac{l}{m},  j+ \frac{l+1}{m} ]
$
\textit{inside} the interval $[m, m+1]$, whence
\begin{eqnarray*}
\sum_{s\in S} W\Bigl( \frac{s}{x\mu} \Bigr)
\:\:
\ge
\:\:
\sum_{s\in S_l} W\Bigl( \frac{s}{x\mu} \Bigr)
&
\ge
&
w_0 \frac{\sigma_S x}{mn}
\\
&=&
w_0 \frac{\sigma_S (x\mu)}{ mn \mu}
         \hspace{20pt}
	 \text{since $\mu>1$ and $n>m$.}
\end{eqnarray*}
Combine these two cases for $j$ and we see that
$$
\sum_{s\in S} W\Bigl( \frac{s}{x} \Bigr) \ge \frac{w_0}{n^2} \sigma_S x.
$$
This completes the proof of the Lemma.
\end{proof}

\vsp

\begin{proof}[Proof of Corollary \ref{cor:slow}]
Without loss of generality we can assume that 
$
f(x) = o(\log\log\log x)
$.
Since $f$ is unbounded and increasing, we can find a sequence
$
0 < x_1 < x_2 < \cdots
$
such that
\begin{equation}
\text{$x_n/ f(x_{n+1}) \rarr 0$ as $n\rarr \infty$.}          \label{xn}
\end{equation}
Define a function $g$ on $\R$ by
$$
g(x)
=
\Bigl\{
  \begin{array}{ll}
      f(x_1) & x< x_2,
      \\
      f(x_i) & x_{i+1} \le x < x_{i+2}, i\ge 1.
  \end{array}
$$
By (\ref{xn}) we have
$
g(x) = o(f(x)) = o(\log\log\log x)
$.
Finally, set
$
k(x) = g(\sqrt{x})
$,
so
\begin{equation}
k(x) = o(f(\sqrt{x})) = o(f(x^{k(x)/2})).          \label{kx}
\end{equation}
On the other hand,
\begin{eqnarray*}
&&
\Bigl( k(x) + \frac{1}{2} + \frac{1}{\sqrt{3}} \Bigr)^{k(x)}
\sum_D
	W\Bigl( \frac{D}{X_{k(x)}} \Bigr)
\\
&
\gg_{E, W}
&
\sum_D
	r_{\text{an}}(E_D)^{k(x)} W\Bigl( \frac{D}{X_{k(x)}} \Bigr)
         \hspace{147pt}
	 \text{by the Main Theorem}
\\
&
\ge
&
f(x^{k(x)/2})^{k(x)}
  \sum_{\substack{D> x^{k(x)/2}\\r_{\text{an}}(E_D)>f(D)}}
W\Bigl(
  \frac{D}{X_{k(x)}}
\Bigr)
         \hspace{97pt}
	 \text{$f$ is increasing}
\\
&
\gg_{E, W}
&
f(x^{k(x)/2})^{k(x)}
X_{k(x)}
\cdot
\Bigl[
  \begin{array}{ll}
     \text{lower density of the set of}
     \\
     \text{$D$ with $r_{\text{an}}(E_D)\ge f(D)$}
  \end{array}
\Bigr],
         \hspace{20pt}
	 \text{by Lemma \ref{lem:lower}.}
\end{eqnarray*}
In light of (\ref{kx}) this lower density must be zero, and Corollary
\ref{cor:slow} follows.
\end{proof}

\vsp

\section{Poisson summation}
      \label{sec:devote}

In this section we adopt Heath-Brown's argument to reduce Proposition
\ref{prop:devoted} to a `multivariable
prime number theorem' for elliptic curves, to be proved in section
\ref{sec:prop2}.

We begin with an auxiliary result.  Denote by $\hat{W_l}$ the Fourier transform
with respect to $t$ of
$$
W_l(x, t, X_k)
:=
\bigl( \log(t^2 X_k^2) + (\log x)/2 \bigr)^l W(t).
$$
Note that the integral defining $\hat{W}_l$ makes sense since $W(0)=0$.  

\begin{lem}
  \label{lem:thrice}
There exists a constant $\gamma_W>0$ depending on $W$ only, so that
for $l>0$, $m\not=0$ and $X_k>2$, as $t\rarr\infty$,
\begin{list}
	{{\rm (\roman{property})}}{\usecounter{property}
				\setlength{\labelwidth}{20pt}}
\item
$|W|, |\hat{W}_l|$ and $|\hat{W}_l|$ all satisfy
$
\displaystyle  < \gamma_W  l^3 (\log X_k + \log x)^l \min(1, |t|^{-3})
$;
\item 
$
\displaystyle
\int_2^{x^r}
  \Bigl|
    \frac{\partial}{\partial t}
    \Bigl( \hat{W}_l\Bigl(x, \frac{X_k m}{t}, X_k\Bigr) \frac{1}{\sqrt{t}} \Bigr)
  \Bigr|
  dt
<
\gamma_W
    l^3 ( \log X_k + \log x )^l
  (T|m|)^{-1/2} \min\Bigl( 1, \Bigl(\frac{x^r}{X_k|m|}\Bigr)^{3/2} \Bigr)
$.
\end{list}
\end{lem}

\begin{proof}
For the rest of this proof, $\gamma_i$ denotes a constant depending on $W$
only.
Since $W(t)=0$ is zero around an open neighborhood of $0$ and since $W$ has
compact support,
$$
\frac{\partial^3}{\partial t^3} W_l(x, t, X_k)
<
\gamma_1 l^3 ( \log X_k + \log x )^l.
$$
Apply 
integration by parts three times and recall that $W$ has compact support, we
get
\begin{eqnarray*}
\hat{W}(x, t, X_k)
&<&
\gamma_2
\frac{1}{|t|^3}
\int_{-\infty}^\infty 
          \frac{\partial^3}{\partial y^3} W_l(x, y, X_k)  dy
\\
&<&
\gamma_3 \, l^3 ( \log X_k + \log x )^l \min(1, |t|^{-3}).
\end{eqnarray*}
The same argument yields the same estimate for
$
\frac{\partial}{\partial t} \hat{W}_l(x, t, X_k)
$
(with different constant).  Consequently,
\begin{eqnarray*}
&&
\frac{\partial}{\partial t}
\Bigl[
  \hat{W}_l\Bigl(x, \frac{X_k m}{t}, X_k\Bigr) \frac{1}{\sqrt{t}}
\Bigr]
\\
&=&
\Bigl(\frac{\partial}{\partial t}\hat{W}_l\Bigr)\Bigl(x, \frac{X_k m}{t}, X_k\Bigr)
\frac{X_k m}{t^{5/2}}
+
\hat{W}_l\Bigl(x, \frac{X_k m}{t}, X_k \Bigr) t^{-3/2}
         \hspace{50pt}
	 \text{chain rule}
\\
&<&
\biggl\{   
\renewcommand{\arraystretch}{1.8}
  \begin{array}{lll}
\displaystyle
    \gamma_4 
    l^3 ( \log X_k + \log x )^l
    \Bigl[
      \Bigl(\frac{X_k m}{t}\Bigr)^{-3} \frac{X_k m}{t^{t/2}}
      +
      \Bigl(\frac{X_k m}{t}\Bigr)^{-3} t^{-3/2}
    \Bigr]
    &&
    \text{if $|X_k m/t|\ge 1$,}
\\
\displaystyle
    \gamma_5
    l^3 ( \log X_k + \log x )^l
    \Bigl[
      \frac{X_k m}{t^{5/2}} + t^{-3/2}
    \Bigr]
    &&
    \text{if $|X_k m/t|<1$}
  \end{array}
\renewcommand{\arraystretch}{1}
\\
&<&
\gamma_6
l^3 ( \log X_k + \log x )^l
t^{-3/2} \min\Bigl( 1, \Bigl|\frac{X_k m}{t}\Bigr|^{-2} \Bigr).
\end{eqnarray*}
So if $|X_k m|\ge x^r$, the integral in the Lemma becomes
$$
< \gamma_7
    l^3 ( \log X_k + \log x )^l
\int_2^{x^r} t^{-3/2} \frac{t^2}{|X_k m|^2} dt
< \gamma_8
    l^3 ( \log X_k + \log x )^l
\frac{x^{3r/2}}{|X_k m|^2}.
$$
On the other hand, if $|X_k m|\le x^r$, then splitting the integral as
$
\int_2^{X_k|m|} + \int_{X_k|m|}^{x^r}
$
gives
$$
< \gamma_9
    l^3 ( \log X_k + \log x )^l
\Bigl( (X_k|m|)^{-1/2} + \int_{X_k|m|}^{x^r} t^{-3/2} dt \Bigr)
< \gamma_{10}
    l^3 ( \log X_k + \log x )^l (X_k|m|)^{-1/2}.
$$
Take $\gamma_W$ to be the maximum of the $\gamma_i$ and the Lemma follows.
\end{proof}

\vsp

Recall the definition of $R(x,D)^r$ and we get
\begin{eqnarray}
&&
\sum_D f(x,D)^i R(x,D)^r W\Bigl(\frac{D}{X_k}\Bigr)
	\nonumber
\\
&=&
2^r \sum_D f(x,D)^i W\Bigl(\frac{D}{X_k}\Bigr)
    \sum_{p_1, \ldots, p_r > 2}
      \beta_{p_1} \cdots \beta_{p_r}
      \Bigl( \frac{D}{p_1} \Bigr)
          \cdots
      \Bigl( \frac{D}{p_r} \Bigr).		\label{number}
\end{eqnarray}

Note that the primes $p_1, \ldots p_r$ in the inner-sum above need not be
distinct.  In particular, the product of the quadratic symbols is a non-trivial
character precisely when $p_1 \cdots p_r$ is not a square.  We proceed
accordingly.

\vsp

\noindent
\framebox{Contribution to (\ref{number}) from those $(p_1, \ldots, p_r)$
whose product is a square}

Then every prime in the $r$-tuple appears with even multiplicity, which means
(i)
$r$ is even, and
(ii)
the product of quadratic characters in (\ref{number}) is
$1$ if every
$
p_i\nmid D
$,
and is zero otherwise.  Thus the contribution in question is
\begin{eqnarray}
&&
2^r
\sum_{p_1, \ldots, p_{r/2}}
      (\beta_{p_1} \cdots \beta_{p_{r/2}})^2
\sum_{D\not\equiv\submod{0}{p_i}}
      f(x,D)^i W\Bigl(\frac{D}{X_k}\Bigr)
      \nonumber
\\
&=&
2^r
\sum_{p_1, \ldots, p_{r/2}}
      (\beta_{p_1} \cdots \beta_{p_{r/2}})^2
\sum_{\delta|\pi'}
      \mu(\delta)
\sum_d
      f(x,d\delta)^i W\Bigl(\frac{d\delta}{X_k}\Bigr),
     \label{square}
\end{eqnarray}
where
$
\pi'= p_1 \cdots p_{r/2}
$
and
$\mu$ = M\"obius function.  The terms in (\ref{square}) with $\delta=1$
sum to
\begin{eqnarray}
&&
2^r
\sum_{p_1, \ldots, p_{r/2}}
      (\beta_{p_1} \cdots \beta_{p_{r/2}})^2
\sum_d
      f(x,d)^i W\Bigl(\frac{d}{X_k}\Bigr)
\nonumber
\\
&=&
2^r \Bigl( \sum_p \beta_p^2 \Bigr)^{r/2}
    \sum_d ( 2\log |d| + (\log x)/2)^i
    W\Bigl( \frac{d}{X_k} \Bigr).
\label{thisis}
\end{eqnarray}
By Lemma \ref{lem:rankin} and Lemma \ref{lem:firstterm}, this is
$$
\le
\:\:
\bigl( 2\log X_k + (\log x)/2 + o(\log x)\bigr)^i
\bigl( 1/3  +o_{E,r}(1)\bigr)^{r/2}
\log^r x
    \sum_d
    W\Bigl( \frac{d}{X_k} \Bigr).
$$
On the other hand, the terms in (\ref{square}) with $\delta>1$ sum to
\begin{eqnarray}
&\ll&
2^r
\sum_{p_1, \ldots, p_{r/2}}
      (\beta_{p_1} \cdots \beta_{p_{r/2}})^2
\sum_{\substack{\delta|\pi'\\ \delta>1}}
      \sum_d
      f(x,d)^i W\Bigl(\frac{d\delta}{X_k}\Bigr)        \nonumber
\\
&\ll&
2^r
\Bigl(2\log X_k + \frac{\log x}{2} + o(1)\Bigr)^i
\!\!\!
\sum_{p_1, \ldots, p_{r/2}}
      (\beta_{p_1} \cdots \beta_{p_{r/2}})^2
\sum_{\substack{\delta|\pi'\\ \delta>1}}
\Bigl[
  \sum_{|d|\le X_k/\delta} 1
  +
  \!\!
  \sum_{|d|>X_k/\delta} \Bigl( \frac{X_k}{d\delta}\Bigr)^{3}
\Bigr]
        \nonumber
\\
&\ll&
2^r
\Bigl(2\log X_k + \frac{\log x}{2} + o(1)\Bigr)^i
\!\!\!
\sum_{p_1, \ldots, p_{r/2}}
      (\beta_{p_1} \cdots \beta_{p_{r/2}})^2
\sum_{\substack{\delta|\pi'\\ \delta>1}}
      X_k/\delta,
                    \label{substack}
\end{eqnarray}
where in the second line we use Lemma \ref{lem:thrice}(a).
The number of $\delta|\pi' = p_1 \cdots p_{r/2}$ is $\le 2^{r/2}$, so 
(\ref{substack}) is
\begin{eqnarray*}
&\ll&
X_k 2^{3r/2} \bigl(2\log X_k + (\log x)/2 + o_{E,i}(1)\bigr)^i
\sum_p \frac{\beta_p^2}{p}
\Bigl(
  \sum_q \beta_q^2
\Bigr)^{r/2-1}
\\
&\ll&
X_k 2^{r/2}
(1/3+ o_E(1))^{r/2-1} \bigl(2\log X_k + (\log x)/2 + o_{E, i}(1) \bigr)^i \log^{r-2} x.
\end{eqnarray*}
Keeping in mind that
$
\sum_D W(D/X_k) \ll_W X_k
$,
we see that if $r$ is even, then the terms in (\ref{number}) coming from those
$
(p_1, \ldots, p_r)
$
whose product is a square, is
$$
\Bigl( 2\log X_k + \frac{\log x}{2} + o_{E, i}(\log x)\Bigr)^i
\bigl( 1/3 +o_E(1)\bigr)^{r/2}
\Bigl(
  \log^r x + O_W(2^{r/2}\log^{r-2} x)
\Bigr)
    \sum_d
    W\Bigl( \frac{d}{X_k} \Bigr).
$$

\vsp

\noindent
\framebox{Contribution to (\ref{number}) from those $(p_1, \ldots, p_r)$
whose product is not a square}

Set
\begin{equation}
\left\{
\begin{array}{llllll}
\pi &=& p_1 \cdots p_r,
\\
\pi_0 &=& \text{largest perfect square divisor of $\pi$ such that
		$(\pi, \pi/\pi_0)=1$,}
\\
\pi_1 &=& \text{the square-free part of $\pi_0$,}
\\
\pi_2 &=& \text{the square-free part of $\pi_1$.}
\end{array}
\right.              \label{pi}
\end{equation}
Then the contribution in question is equal to
$$
2^r \sum_{\substack{p_1, \ldots, p_r\\ \pi_2>1}}
      \beta_{p_1} \cdots \beta_{p_r}
    \sum_{\submod{j}{\pi_1 \pi_2}}
	\Bigl( \frac{j}{\pi_2} \Bigr)
    \sum_{m=-\infty}^\infty
	f(x, j+m\pi_1\pi_2)^i W\Bigl( \frac{j+m\pi_1 \pi_2}{X_k} \Bigr).
$$
Set
$
e(z)=\exp(2\pi i z)
$.
Apply Poisson summation and we get 
\begin{eqnarray}
&&
2^r \sum_{\substack{p_1, \ldots, p_r\\ \pi_2>1}}
      \beta_{p_1} \cdots \beta_{p_r}
    \sum_{\submod{j}{\pi_1 \pi_2}}
	\Bigl( \frac{j}{\pi_2} \Bigr)
    \sum_{m=-\infty}^\infty
	\hat{W}_i\Bigl(x, \frac{X_k m}{\pi_1\pi_2}, X_k\Bigr)
	\frac{X_k}{\pi_1 \pi_2}
	e\Bigl( \frac{mj}{\pi_1 \pi_2} \Bigr)
\nonumber
\\
&=&
2^r X_k \sum_{\substack{p_1, \ldots, p_r\\ \pi_2>1}}
      \frac{\beta_{p_1} \cdots \beta_{p_r}}{\pi_1 \pi_2}
    \sum_{m=-\infty}^\infty
	\hat{W}_i\Bigl(x, \frac{X_k m}{\pi_1\pi_2}, X_k\Bigr)
    \sum_{\submod{j}{\pi_1 \pi_2}}
	\Bigl( \frac{j}{\pi_2} \Bigr)
	e\Bigl( \frac{mj}{\pi_1 \pi_2} \Bigr).
\label{lineone}
\end{eqnarray}
Since $\pi_2>1$, if $\pi_1\pi_2| m$ then the $j$-sum in (\ref{lineone})
is zero.  So suppose
$
\pi_1\pi_2\nmid m
$;
in particular, $m\not=0$.  For $l=1, 2$, set
$$
\delta_l = (\pi_l, m),
\:\:
\pi_l = \delta_l \pi_l',
\:\:
m = \delta_l n_l.
$$
Since $(\pi_1, \pi_2)=1$, by the Chinese remainder theorem the $j$-sum in
(\ref{lineone}) is
\begin{eqnarray}
&&
\Bigl[
\sum_{\submod{j_1}{\pi_1}}
	e\Bigl( \frac{mj_1}{\pi_1} \Bigr)
\Bigr]
\Bigl[
\sum_{\submod{j_2}{\pi_2}}
	\Bigl( \frac{j_2}{\pi_2} \Bigr)
	e\Bigl( \frac{mj_2}{\pi_2} \Bigr)
\Bigr]
\nonumber
\\
&=&
\Bigl[
\sum_{\submod{l_1}{\pi_1'}}
	e\Bigl( \frac{n_1 l_1}{\pi_1'} \Bigr)
\sum_{\substack{\submod{j_1}{\pi_1}\\ j_1\equiv\submod{l_1}{\pi_1'}}}
	1
\Bigr]
\Bigl[
\sum_{\submod{l_2}{\pi_2'}}
	\Bigl( \frac{l_2}{\pi_2'} \Bigr)
	e\Bigl( \frac{n_2 l_2}{\pi_2'} \Bigr)
\sum_{\substack{\submod{j_2}{\pi_2}\\j_2\equiv\submod{l_2}{\pi_2'}}}
	\Bigl( \frac{j_2}{\delta_2} \Bigr)
\Bigr].
		\label{j2}
\end{eqnarray}
Note that the $j_2$-sum in (\ref{j2}) is zero unless $\delta_2=1$, and the
$j_1$-sum is
$
\delta_1
$.
Moreover, $\pi_1$, and hence
$
\pi_1'
$,
is square-free, so (\ref{j2}) is
\begin{eqnarray*}
&=&
(-1)^{\#\{p|\pi_1'\}} \delta_1
\sum_{\submod{j}{\pi_2}}
	\Bigl( \frac{j}{\pi_2} \Bigr)
	e\Bigl( \frac{n\delta_1 j}{\pi_2} \Bigr)
\\
&=&
(-1)^{\#\{p|\pi_1'\}} 
\frac{\delta_1\sqrt{\pi_2}}{1+\sqrt{-1}}
\Bigl(\frac{n\delta_1}{\pi_2}\Bigr)
\Bigl(
  1 - \sqrt{-1} \Bigl( \frac{-1}{\pi_2} \Bigr)
\Bigr),
\end{eqnarray*}
by the standard quadratic Gauss sum calculation.  Recall that 
$
m=n\delta_1\not=0
$
and we see that (\ref{lineone}) is
\begin{eqnarray}
&\ll&
2^r X_k
\sum_{\substack{p_1, \ldots, p_r\\ \pi_2>1}}
	\frac{\beta_{p_1}\cdots \beta_{p_r}}{\pi_1\sqrt{\pi_2}}
\sum_{\delta_1|\pi_1}
\sum_{\substack{|n|\not=0\\p_j\nmid n}}
	\hat{W}_i
		\Bigl(
		x, \frac{T\delta_1 n}{\pi_1 \pi_2}, X_k
		\Bigr)
\delta_1
	\Bigl(
		\frac{\pm n\delta_1}{\pi_2}
	\Bigr)
\nonumber
\\
&\ll&
2^r X_k
\sum_{|n|\not=0}
\Bigl|
	\sum_{\substack{p_1, \ldots, p_r\\ \pi_2>1\\(p_j, n)=1}}
		\beta_{p_1}\cdots \beta_{p_r}
	\sum_{\delta_1|\pi_1}
	\Bigl(
		\frac{\pm n\delta_1}{\pi_2}
	\Bigr)
	\frac{1}{\pi_1' \sqrt{\pi_2}}
	\hat{W}_i
		\Bigl(
		x, \frac{Tn}{\pi_1' \pi_2}, X_k
		\Bigr)
\Bigr|.
          \label{trivial}
\end{eqnarray}
We now  estimate (\ref{trivial}) in two ways, first unconditionally
and then invoke the GRH.

\vsp

\noindent
\framebox{Unconditional Estimate}

For the unconditional estimate we will take the test function $F$ to be $F_3$,
in which case
$
|| F_3 || \le 1
$,
whence
$
|\beta_p| \le (2\log p)/\sqrt{p}
$.
There are
$
\le 2^r
$
terms in the $\delta$-sum in
$
Q(p_1, \ldots, p_r, n)
$.
Since $F_3$ vanishes outside $(-1,1)$, we have $\beta_p=0$ if
$p>x$.  Use Lemma \ref{lem:thrice}(a) to bound
$
\hat{W}_i
$
and we see that
\begin{eqnarray*}
(\ref{trivial})
&
\ll_W
&
2^r X_k \sum_{|n|\not=0} 
  \sum_{p_1, \ldots, p_r <x}
    2^r
    \frac{\log p_1\cdots \log p_r}{p_1 \cdots p_r}
  \frac{(p_1 \cdots p_r)^3}{T^3 |n|^3}
  i^3 (\log X_k + \log x)^i
\\
&\ll_W&
4^r i^3 (\log X_k + \log x)^{r+i} x^{3r}/X_k^2.
\end{eqnarray*}

\vsp

\noindent
\framebox{GRH Estimate}

First, rewrite (\ref{trivial}) as 
\begin{equation}
2^r X_k
\sum_{|n|\not=0}
\Bigl|
    \sum_{u\ge 2}
        \hat{W}_i
		\Bigl(
		x, \frac{Tn}{u}, X_k
		\Bigr)
        \frac{1}{\sqrt{u}}
	\sum_{\substack{p_1, \ldots, p_r\\ \pi_2>1\\(p_j, n)=1\\ p_1 \cdots p_r=u}}
	\underbrace{\beta_{p_1}\cdots \beta_{p_r}
	\sum_{\delta_1|\pi_1}
	\Bigl(
		\frac{\pm n\delta_1}{\pi_2}
	\Bigr)
	\frac{1}{\sqrt{\pi_1'}}}_{:= Q(p_1, \ldots, p_r, n)}
\Bigr|
         \label{absvalue}
\end{equation}
Note that if
$
p_1 \cdots p_r \ge x^r
$,
then $Q(p_1, \ldots, p_r, n) = 0$  for any $n$.  In particular, the $u$-sum
in (\ref{absvalue}) is a finite sum.  To evaluate this $u$-sum we proceed
by partial summation.  
That calls for the following estimate, to be proved in sections
\ref{sec:key} and \ref{sec:prop2}.

\vsp

\begin{prop}
             \label{prop:fundamental}
Assume the GRH for every $L(E_D, s)$.   Then there exists a constant
$\tilde{c}_E$ depending on $E$ only so that, for any integers $m, r>0$, as
$
p_1, \ldots, p_r
$
run through all prime numbers,
\begin{eqnarray}
\sum_{\substack{p_1 \cdots p_r \le U\\ \pi_2>1\\(p_j,n)=1}}
	Q(p_1, \ldots, p_r, n)
&
\ll
&
\tilde{c}_E^r \Bigl[\log N_E + 3\log |n| + 3\log(U+2)\Bigr]^r \log^{2r+1}x.
\end{eqnarray}
\end{prop}

\vsp

\noindent
Assuming this, the $u$-sum in (\ref{absvalue}) is
\begin{eqnarray}
&=&
\Bigl[
  \sum_{\substack{p_1 \cdots p_r\le x^r\\(p_j,n)=1\\ \pi_2>1}}
  Q(p_1, \ldots, p_r, n)
\Bigr]
\hat{W}_i\Bigl( x, \frac{X_k n}{x^r}, X_k \Bigr)
\frac{1}{\sqrt{x^r}}
                       \label{partial1}
\\
&&
-
\int_2^{x^r}
\Bigl[
  \sum_{\substack{p_1 \cdots p_r\le t\\(p_j,n)=1\\ \pi_2>1}}
  Q(p_1, \ldots, p_r, n)
\Bigr]
\frac{\partial}{\partial t}
\Bigl(
  \hat{W}_i\Bigl( x, \frac{X_k n}{t}, X_k \Bigr) \frac{1}{\sqrt{t}}
\Bigr)
dt
                       \label{partial2}
\\
&\ll_W&
\tilde{c}_E^{\:r} \Bigl[\log N_E + 3\log |n| + 3\log(U+2)\Bigr]^r \log^{2r+1}x
\times
i^3 (\log X_k + \log x)^i
                       \nonumber
\\
&&
\times
\Bigl[
  \frac{1}{\sqrt{x^r}}
  \min\Bigl(1, \Bigl|\frac{x^r}{X_k m}\Bigr|^3\Bigr)
  +
  \frac{1}{\sqrt{X_k |m|}}
  \min\Bigl(1, \Bigl|\frac{x^r}{X_k m}\Bigr|^{\frac{3}{2}}\Bigr)
\Bigr]
                       \nonumber
\\
&
\ll_{E, W}
&
r^3 (3\tilde{c}_E)^r
(\log X_k + \log x)^i \Bigl[\log |n| + \log(U+2)\Bigr]^r  
\frac{\log^{2r+1}x}
     {\sqrt{X_k |m|}}
  \min\Bigl(1, \Bigl|\frac{x^r}{X_k m}\Bigr|^{\frac{3}{2}}\Bigr).
                       \nonumber
\end{eqnarray}
Recall that $U\le x^r$.  Consequently,  (\ref{absvalue}) becomes
\begin{eqnarray}
&
\ll_{E, W}
&
r^{r+3} (3\tilde{c}_E)^r
(\log X_k + \log x)^i  (\log |n|  + \log x)^r
\sum_{|n|\not=0}
	\frac{\sqrt{X_k}}{\sqrt{|m|}}
	\min\Bigl( 1, \Bigl|\frac{x}{X_k m}\Bigr|^{\frac{3}{2}} \Bigr).
                            \label{become}
\end{eqnarray}
Thus the contribution to the $n$-sum from those $|n|\ge x^r/X_k$ is 
\begin{eqnarray*}
&
\ll_{E,  W}
&
r^{r+3} (3\tilde{c}_E)^r
(\log X_k + \log x)^i  (\log |n|  + \log x)^r
\sqrt{X_k}\sum_{|n|\ge x^r/X_k} \frac{1}{\sqrt{n}}\Bigl( \frac{x}{X_k m} \Bigr)^{3/2}
\\
&
\ll_{E, W}
&
r^{r+3} (3\tilde{c}_E)^r
(\log X_k + \log x)^{r+i}
\sqrt{X_k} \Bigl(\frac{x}{X_k}\Bigr)^{3/2} \sum_{|n|\ge x/X_k} \frac{\log^r |n|}{n^2}
\\
&
\ll_{E, W}
&
r^{r+3} (3\tilde{c}_E)^r
(\log X_k + \log x)^{r+i} x^{r/2}.
\end{eqnarray*}
On the other hand, the contribution from those $|n|<x^r/X_k$ is
\begin{eqnarray*}
&
\ll_{E, W}
&
r^{r+3} (3\tilde{c}_E)^r
(\log X_k + \log x)^{r+i}
\sqrt{X_k}
        \sum_{0<|n|<x/X_k}
        \frac{\log^r |n|}{\sqrt{|n|}}
\\
&
\ll_{E,  W}
&
r^{r+3} (3\tilde{c}_E)^r
(\log X_k + \log x)^{r+i} x^{r/2}.
\end{eqnarray*}
This completes the proof of Proposition \ref{prop:devoted}.
\qed

\vsp

\begin{remark}
      \label{remark:cubic}
The argument in this section readily extends to higher order twist families.
The main difference, say for the cubic twist family
$
\mathcal{E}_m: x^3 + y^3 = m
$,
is that the argument now proceeds according to whether 
$
p_1 \cdots p_r
$
is a perfect cube or not.  The rest of the argument, including Proposition
\ref{prop:fundamental}, extends with no change.  As a result, the Main Theorem
extends to the cubic twist family $\mathcal E_m$
with the factor $1/2$ replaced by $1/3$.
\end{remark}

\begin{remark}
      \label{remark:failure}
While Proposition \ref{prop:fundamental} gives an essentially optimal bound
for the size of the $Q$-sum, we have no control over the \textit{sign} of
this $Q$-sum as $u$ varies.  Because of that, to estimate (\ref{partial1}) and
(\ref{partial2}) using Proposition \ref{prop:fundamental} we are forced to put
absolute value signs everywhere.  This is essentially the only place in the
proof of the Main Theorem where we might lose information
(the $\ll$ in (\ref{trivial}) does not have any material impact on the rest of
the proof).
\end{remark}

\vsp

\section{A complex prime number theorem}
    \label{sec:key}

The results in this section are elliptic curves analog of classical estimates;
we provide the details for lack of a good reference.
As is customary, given a complex number $s$ we denote by $\sigma$ and $t$
its real and imaginary part, respectively.

\begin{lem}
  \label{lem:log}
Assume the GRH for $L(E,s)$.  Then for $\sigma\ge 1 + 1/\log x$ and $|t|\ge 2$,
we have the estimate
$$
L'(E,s)/L(E,s) \ll (\log N_E + \log (|s|+2) ) \log x.
$$
\end{lem}

\begin{proof}
We have the basic relation
\begin{equation}
-\frac{L'(E, s)}{L(E,s)}
=
\log\frac{\sqrt{N_E}}{2\pi}
+
\frac{\Gamma'(s)}{\Gamma(s)}
-
B_E
-
\sum_\rho \Bigl( \frac{1}{s-\rho} + \frac{1}{\rho} \Bigr),
      \label{log}
\end{equation}
where $B_E$ is a constant depending only on $E$, and $\rho$ runs through the
non-trivial zeros of $L(E, s)$.  Since
$
\ov{L(E, s)} = L(E,\ov{s})
$,
complex conjugation takes the zeros of $L(E,s)$ to themselves; from (\ref{log})
we see that $B_E$ is real, and that as in \cite[p.~83]{davenport},
$$
B_E = -\sum_\rho Re\Bigl(\frac{1}{\rho} \Bigr).
$$
The $\Gamma$-term in (\ref{log}) is
$
\ll \log |t|
$
if $|t|\ge 2$ and $1\le \sigma \le 3$.  It follows that
$$
Re\Bigl(-\frac{L'(E, s)}{L(E,s)}\Bigr)
\ll
\Bigl( \log N_E + \log (|t|+2) \Bigr)
-
\sum_\rho Re\Bigl( \frac{1}{s-\rho} \Bigr).
$$
Since $L'(E,s)/L(E,s)$ is bounded on the line Re$(s)=2$, for such $s$ we get
$$
\sum_\rho \text{Re}\Bigl( \frac{1}{s-\rho} \Bigr)
\ll
\log N_E + \log (|t|+2).
$$
Write $\rho = \beta + i\gamma$.  Then for $s=2+it$,
$$
\text{Re}\frac{1}{s-\rho}
=
\frac{2-\beta}{(2-\beta)^2 + (t-\gamma)^2}
\ge
\frac{1/2}{1 + (t-\gamma)^2},
$$
whence
\begin{equation}
\sum_\rho \frac{1}{1+(t-\gamma)^2}
\ll
\log N_E + \log(|t|+2).
    \label{sumrho}
\end{equation}
Standard argument then shows that 
$$
\frac{L'(E,s)}{L(E,s)} = \sideset{}{'}\sum_\rho \frac{1}{s-\rho}
+
O(\log N_E + \log(|t|+2)),
$$
where the sum runs over those $\rho$ for which $|T-\gamma|<1$.  By
(\ref{sumrho}) there are
$
\ll \log N_E + \log(|t|+2)
$
such $\rho$, and under GRH,
$
|s-\rho| \ge 1/\log x
$
if $\sigma\ge 1+1/\log x$.  The Lemma then follows.
\end{proof}

\vsp

\begin{lem}
  \label{lem:basic}
Assume the GRH for $L(E, s)$.  For $j\ge 0, x\gg_E 1$ and 
$
1 + 1/\log x \le \sigma \le 2
$,
we have the estimate
\begin{eqnarray*}
\frac{1}{\log^{j} x} \sum_{p<x} \frac{a_p(E)\log^{1+j} p}{p^s}    
\ll
(\log N_E + \log(|s|+2) ) \log^2 x.
\end{eqnarray*}
\end{lem}

\vsp

Recall the definition of $F$ and we get immediately

\vsp

\begin{cor}
  \label{cor:keep}
Assume the GRH for $L(E, s)$.  Then there exists a constant $\epsilon_E>0$ such that,
for
$
1 + 1/\log x \le \sigma \le 2
$,
we have the estimate
\begin{eqnarray*}
\sum_{p<x}
    \frac{a_p(E)\log p}{p^s}  F\Bigl(\frac{\log p}{\log x}\Bigr)
\le
\epsilon_E (\log N_E + \log(|s|+2) ) \log^2 x.           \qed
\end{eqnarray*}
\end{cor}

\begin{proof}[Proof of Lemma \ref{lem:basic}]
By partial summation it suffices to take $j=0$.  To handle that case we mimic
the proof of the prime number theorem under the Riemann hypothesis.

Recall the definition of 
$
c_n(E)
$
in section \ref{sec:quad}.
Set $c=1/2+1/\log x$.  Apply the Perron formula \cite[(2) on p.~104]{davenport}
and we get, for $\sigma > 1$,
\begin{eqnarray}
\lefteqn{\Bigl|
  \int_{c-\sqrt{-x}}^{c+\sqrt{-x}}
  \frac{L'(E,\sigma+it+\xi)}{L(E,\sigma+it+\xi)} \frac{x^\xi}{\xi}d\xi
  -
  \sum_{n<x}^\infty \frac{c_n(E)\log n}{n^{\sigma+it}}
\Bigr|}
           \nonumber
\\
& \ll &
\sum_{\substack{n=1\\n\not=x}}
        \frac{\Lambda(n)}{n^{\sigma-1/2}}
	\Bigl(
	  \frac{x}{n}
	\Bigr)^c
	\min\Bigl(1, \frac{1}{\sqrt{x} \bigl|\log \frac{x}{n}\bigr|}\Bigr)
	+
	\frac{c\Lambda(n)}{\sqrt{x} n^{\sigma-1/2}},
\label{perron}
\end{eqnarray}
where $\Lambda$ denotes the usual von Mongoldt function, and the last term on
the right side of (\ref{perron}) is present only if $x$ is a prime power.

If $n \ge \frac{5}{4}x$ or if $n\le \frac{3}{4}x$ then
$
|\log \frac{x}{n}|
$
has a positive lower bound.  Thus the contribution of such $n$ to the right
side of (\ref{perron}) is (recall that $\sigma> 1$)
$$
\ll
\sum_n \frac{\Lambda(n)}{n^{1+1/\log x}}
\ll
\Bigl(\frac{-\zeta'(1+1/\log x)}{\zeta(1+1/\log x)}\Bigr)
\ll
\log x.
$$
The argument in \cite[p.~107]{davenport} shows that the contribution from
those $n$ such that
$
\frac{3}{4}x < n < \frac{5}{4}x
$,
$x\not=$ prime power, is
$$
\ll
\frac{\log x}{\sqrt{x}} \min\Bigl(1, \frac{x}{\sqrt{x}\langle x\rangle}\Bigr)
+
\log^2 x.
$$
Putting everything together and we get
$$
\Bigl|
  \int_{c-\sqrt{-x}}^{c+\sqrt{-x}}
  \frac{L'(E,\sigma+it+\xi)}{L(E,\sigma+it+\xi)} \frac{x^\xi}{\xi}d\xi
  -
  \sum_{n<x} \frac{c_n(E)\log n}{n^{\sigma+it}}
\Bigr|
\ll
\log^2 x
+
\frac{\log x}{\sqrt{x}} \min\Bigl(1, \frac{x}{\sqrt{x}\langle x\rangle}\Bigr).
$$
Our next step is to estimate the integral.  Since $\sigma \ge 1$, under the
GRH the integrand has no pole inside the rectangle with vertices
$$
c+\sigma+it \pm iT,  1 + \frac{1}{\log x} + it \pm i\sqrt{x}.
$$
Thus it remains to estimate the integral along the other three edges of
this rectangle.

The integral along the top edge is (recall $1<\sigma\le 2$)
\begin{eqnarray*}
&&
\int_{c}^{-\sigma+1/\log x}
  \frac{-L'(E, \sigma+it+\xi+i\sqrt{x})}{L(E, \sigma+it+\xi+i\sqrt{x})}
  \frac{x^{\xi+i\sqrt{x}}}{\xi+i\sqrt{x}} d\xi
\\
&\ll&
\frac{\sqrt{x}}{\sqrt{x}}
\int_{c}^{-\sigma+1/\log x}
  \bigl[
    \log N_E + \log( |\xi + \sigma + it + i\sqrt{x}| + 2)
  \bigr]
  \log x
\:
  d\xi
        \hspace{20pt}
	\text{by Lemma \ref{lem:log}}
\\
&\ll&
\frac{\sqrt{x}\log x}{\sqrt{x}}
  \bigl[
    \log N_E + \log(|t|+2)
  \bigr].
\end{eqnarray*}
The same bound holds for the integral along the bottom edge.  As for the
vertical edge,
\begin{eqnarray*}
&&
\int_{-\sqrt{x}}^{\sqrt{x}}
  \frac{L'(E, 1+\frac{1}{\log x}+it+i\tau)}{L(E, 1+\frac{1}{\log x}+it+i\tau)}
  \frac{x^{1/\log x}}{\frac{1}{\log x} + i\tau}
  d\tau
\\
&\ll&
  \bigl(
    \log N_E + \log( 1+\frac{1}{\log x} + |t|+\sqrt{x}+2)
  \bigr)
  \log x
  \Bigl(
    \int_0^2 \log x \: d\tau
    +
    \int_2^{\sqrt{x}} \frac{d\tau}{\tau}
  \Bigr)
\\
&\ll&
  \bigl(
    \log N_E + \log(|t|+\sqrt{x} +2)
  \bigr)
  \bigl(
    \log^2 x + \log x \log \sqrt{x}
  \bigr).
\end{eqnarray*}
Putting everything together, we get, for $\sigma\ge 1 + 1/\log x$,
$$
\sum_{n<x} \frac{c_n(E)\log n}{n^{\sigma+it}}
\ll
\log^2 x
+
\frac{\log x}{\sqrt{x}}
+
\bigl(
  \log N_E + \log(|t|+\sqrt{x}+2)
\bigr)
\log^2 x.
$$
Since $\sigma>1$, the contribution to the sum on the left side from non-prime
$n$ is
$
\ll \sum_{m<\sqrt{x}} \frac{\log m}{m^{3/2}} \ll 1
$,
so we are done.
\end{proof}

\vsp

\section{Proof of Proposition \ref{prop:fundamental}}
     \label{sec:prop2}

When $r=1$, Brumer \cite[(2.13)]{brumer}
 deduces Proposition \ref{prop:fundamental} from the explicit
formula in conjunction with an estimate of a weighted sum of zeros of
$L(E_D, s)$.   Another (essentially equivalent) way is to apply the Perron
formula as in the proof of the prime number theorem to the
logarithmic derivative of $L(E, s)$.  The explicit formula approach does not
seem to generalize to $r>1$, but the approach via the Perron formula does,
with the key analytic estimate provided by Corollary \ref{cor:keep}.  We
prove Proposition  \ref{prop:fundamental} in several steps.

\vsp

\noindent
\framebox{\bf Step I.}   
\hspace{5pt}
Define
\begin{eqnarray*}
L_x(E,s)
=
\sum_{p<x}
    \frac{a_p(E)\log p}{p^{s}} F\Bigl( \frac{\log p}{\log x}\Bigr).
\end{eqnarray*}
This is a finite sum, and hence it is holomorphic for all $s$.  Apply
the Perron formula as in the proof of Lemma \ref{lem:basic}, we get
\begin{eqnarray*}
\biggl|
\int_{\frac{1}{\log x} - \sqrt{-x}}^{\frac{1}{\log x} + \sqrt{-x}}
  L_x(E,s+1)^r
\frac{U^s}{s} ds
-\sum_{p_1 \cdots p_r \le U}  \beta_{p_1}(E) \cdots \beta_{p_r}(E)
\biggr|
&\ll&
\log^2 x.
\end{eqnarray*}
As for the integral, Corollary \ref{cor:keep} shows that it is
$$
\le
\epsilon_E^{\:r}
\bigl(
   \log N_E + \log (U+2)
\bigr)^r
\log^{2r}\!x \:\: U^{1/\log x}
\Bigl[
  \int_0^2 \frac{dt}{\bigl|\frac{1}{\log x}+it\bigr|}
  +
  \int_2^{\sqrt{x}} \frac{dt}{t}
\Bigr].
$$
Recall that $U \le x^r$ and we get
\begin{eqnarray}
\sum_{p_1 \cdots p_r \le U}  \beta_{p_1}(E) \cdots \beta_{p_r}(E)
&\ll&
(\epsilon_E e)^r
\bigl(
   \log N_E + \log (U+2)
\bigr)^r
\log^{2r+1}\!x.                  \label{step1}
\end{eqnarray}

\vsp

\noindent
\framebox{\bf Step II.}
\hspace{5pt}
Fix an integer $m\not=0$.  With $\pi_2$ defined as in (\ref{pi}), we claim that
\begin{eqnarray}
\sum_{\substack{p_1 \cdots p_r \le U\\ \pi_2>1}}
      \beta_{p_1} \cdots \beta_{p_r}
\Bigl(
  \frac{m}{p_1 \cdots p_r}
\Bigr)
&\ll&
(2\epsilon_E e)^r
\bigl(
   \log N_E + 2\log |m| + \log (U+2)
\bigr)^r
\log^{2r+1}\!x.                  \label{step2}
\end{eqnarray}
To say that $\pi_2=1$ means that $r$ is even and
$
\pi = (p_1 \cdots p_{r/2})^2
$,
so
\begin{eqnarray*}
\sum_{\substack{p_1 \cdots p_r \le U\\ \pi_2>1}}
      \beta_{p_1} \cdots \beta_{p_r}
\Bigl(
  \frac{m}{p_1 \cdots p_r}
\Bigr)
&\le&
\sum_{p_1 \cdots p_{r/2} \ll \sqrt{U}}
      \Bigl(
        \beta_{p_1} \cdots \beta_{p_{r/2}}
      \Bigr)^2
\\
&\le&
\Bigl(
  \sum_{p<x} \beta_p^2
\Bigr)^{r/2}
      \hspace{20pt}
      \text{since $\beta_p=0$ if $p\ge x$}
\\
&\ll&
(4\log^2 x)^{r/2},
\end{eqnarray*}
which is satisfactory.  Thus it remains to study (\ref{step2}) without the
additional condition $\pi_2>1$.
If $p\nmid 2N_E m$ then
$
a_p(E) \bigl( \frac{m}{p} \bigr) = a_p(E_m)
$,
so the left side of (\ref{step2}) without the $\pi_2$ condition is
\begin{eqnarray*}
&=&
\sum_{p_1 \cdots p_r \le U}
  \beta_{p_1}(E_m) \cdots \beta_{p_r}(E_m) 
+
\\
&&
\hspace{20pt}
O\Bigl[
  \sum_{j=1}^r
  \sum_{\substack{p_1,  \ldots, p_j\\p_l|2N_E m}}
	\frac{\log p_1 \cdots \log p_j}{\sqrt{p_1 \cdots p_j}}
\hspace{-10pt}
  \sum_{q_1 \cdots q_{r-j} \le U/p_1 \cdots p_j}
\hspace{-10pt}
    \beta_{q_1}(E) \cdots \beta_{q_{r-j}}(E) 
  \Bigl(
	\frac{m}{q_1 \cdots q_{r-j}}
  \Bigr)
\Bigr].
\end{eqnarray*}
We estimate the first sum above using (\ref{step1}), and we estimate each of
the inner $q$-sum in the $O$-term by induction.  All together, this yields
\begin{eqnarray*}
&\ll&
(\epsilon_E e)^r
\bigl[
  \log N_{E_m}  + \log(U+2)
\bigr]^r
\log^{2r+1}x
\\
&&
\hspace{20pt}
+
\sum_{j=1}^r
  2^j
  (\epsilon_E e)^{r-j}
        \Bigl[
	  \log N_{E} + 2\log |m| + \log(U+2)
	\Bigr]^{r-j}
	\log^{2(r-j)+1}\!x
  \Bigl(
        \sum_{p| 2N_E m} \frac{\log p}{\sqrt{p}}
  \Bigr)^j.
\end{eqnarray*}
Back in section \ref{sec:quad} we saw that the $p$-sum is
$
\ll \log^{3/4} (2N_E m)
$.
Also, $N_{E_m} \ll N_E m^2$, and Step II follows.

\vsp

\noindent
\framebox{\bf Step III.}
\hspace{5pt}
Fix an integer $n\not=0, 1$.  We claim that
\begin{eqnarray}
\lefteqn{\sum_{\substack{p_1 \cdots p_r \le U\\(n, p_j)=1\\ \pi_2>1}}
          \beta_{p_1} \cdots \beta_{p_r}
\Bigl(
  \frac{m}{p_1 \cdots p_r}
\Bigr)}                  \label{step3}
\\
&\ll&
(2\epsilon_E e)^r
\bigl(
   \log N_E + 2\log |m| + \log |n| + \log (U+2)
\bigr)^r
\log^{2r+1}\!x.
               \nonumber
\end{eqnarray}

By (\ref{step2}),
\begin{eqnarray*}
\sum_{\substack{p\le U\\(n,p)=1}} \beta_p \Bigl(\frac{m}{p}\Bigr)
&\ll&
(2\epsilon_E e)
\bigl(
   \log N_E + 2\log |m| + \log (U+2)
\bigr)
\log^3\!x
+
\sum_{p|n} \frac{\log p}{\sqrt{p}}
\\
&\ll&
(2\epsilon_E e)
\bigl(
   \log N_E + 2\log |m|  + \log (U+2)
\bigr)
\log^3\!x
+
\log |n|.
\end{eqnarray*}
This gives the case $r=1$.  In general,
\begin{eqnarray*}
\lefteqn{\sum_{\substack{p_1 \cdots p_r \le U\\(n, p_j)=1\\ \pi_2>1}}
          \beta_{p_1} \cdots \beta_{p_r}
\Bigl(
  \frac{m}{p_1 \cdots p_r}
\Bigr)
\:\:
=
\sum_{\substack{p_1 \cdots p_r \le U\\ \pi_2>1}}
          \beta_{p_1} \cdots \beta_{p_r}
\Bigl(
  \frac{m}{p_1 \cdots p_r}
\Bigr)}
\\
&&
+
O\Bigl(
  \sum_{j=1}^r
    \:\:
    2^j
    \hspace{-8pt}
  \sum_{\substack{p_1 \cdots p_j\le U\\p_l|n}}
    \frac{\log p_1\cdots \log p_j}
	 {\sqrt{p_1 \cdots p_j}}
  \hspace{12pt}
  \Bigl|
    \hspace{-20pt}
    \sum_{\substack{q_1\cdots q_{r-1}\le U/p_1\cdots p_j\\(q_l,n)=1\\ \pi_2>1}}
          \beta_{q_1} \cdots \beta_{q_{r-j}}
	  \Bigl(
	    \frac{m}{q_1 \cdots q_{r-j}}
	  \Bigr)
  \Bigr|
  \:
\Bigr).
\end{eqnarray*}
Step III now follows from (\ref{step2}) plus induction on $r$.

\vsp

\noindent
\framebox{\bf Step IV.}
\hspace{5pt}
Finally we come to prove Proposition \ref{prop:fundamental}.   We proceed by
induction on $r$, the case $r=1$ being automatic.

By (\ref{step3}), the sum of terms with $\pi_1=1$ is
\begin{equation}
\ll
(2\epsilon_E e)^r
\bigl(
   \log N_E + 2\log |m| + \log |n| + \log (U+2)
\bigr)^r
\log^{2r+1}\!x.          \label{pione}
\end{equation}
It remains to account for terms with $\pi_1>1$.  That happens precisely when
$\pi$ is exactly divisible by an even prime power.  
Then the contribution from these terms is therefore equal to
($
\lfloor z \rfloor :=
$
the largest integer $\le z$)
\begin{equation}
\sum_{\lambda=1}^{\lfloor r/2\rfloor}
\sum_{\substack{p\le U\\p\nmid n}}
     \beta_p^{2\lambda}
\sum_{\substack{q_1 \cdots q_{r-2\lambda}\le U/p^{2\lambda}\\
                (q_j, np)=1\\
                \pi_2>1}}
     \beta_{q_1} \cdots \beta_{q_r-2\lambda}
\sum_{\delta_1|\pi_1}
     \Bigl(
       \frac{\pm n\delta_1}{\pi_2}
     \Bigr)
     \frac{1}{\sqrt{\pi_1'}},         \label{smile}
\end{equation}
where $\pi_1$ and $\pi_2$ above are defined with respect to the $r$-tuple
$
(q_1, \ldots, q_{r-2\lambda}, \overbrace{p, \ldots, p}^{2\lambda})
$.
If we denote by $\pi_1(q)$ and $\pi_2(q)$ the corresponding quantities in
(\ref{pi}) with respect to the $(r-2\lambda)$-tuple
$
(q_1, \ldots, q_{r-2\lambda})
$,
then  $\pi_2 = \pi_2(q)$  and  $\pi_1 = \pi_1(q)p$, so  (\ref{smile}) is equal
to
\begin{eqnarray}
&&
\sum_{\lambda=1}^{\lfloor r/2\rfloor}
\sum_{\substack{p\le U\\p\nmid n}}
     \beta_p^{2\lambda}
\sum_{\substack{q_1 \cdots q_{r-2\lambda}\le U/p^{2\lambda}\\
                (q_j, np)=1\\
                \pi_2>1}}
     \beta_{q_1} \cdots \beta_{q_r-2\lambda}
\sum_{\delta_1|\pi_1(q)}
     \Bigl[
       \Bigl(
         \frac{\pm n\delta_1}{\pi_2}
       \Bigr)
       \frac{1}{\sqrt{\pi_1'p}}
       +
       \Bigl(
         \frac{\pm n\delta_1 p}{\pi_2}
       \Bigr)
       \frac{1}{\sqrt{\pi_1'}}
     \Bigr]
                    \nonumber
\\
&\ll&
\sum_{\lambda=1}^{\lfloor r/2\rfloor}
\sum_{\substack{p\le x\\p\nmid n}}
     \Bigl(
       \frac{4\log^2 p}{p}
     \Bigr)^\lambda
\Bigl[
     \frac{1}{\sqrt{p}}
     \sum_{q_1 \cdots q_{r-2\lambda}\le U/p^{2\lambda}}
            \beta_{q_1} \cdots \beta_{q_{r-2\lambda}}
     \sum_{\delta_1|\pi_1(q)}
            \Bigl(
	      \frac{\pm n\delta_1}{\pi_2}
	    \Bigr)
       \frac{1}{\sqrt{\pi_1'}}
\label{here}
\\
&&
\hspace{108pt}
     +
     \sum_{q_1 \cdots q_{r-2\lambda}\le U/p^{2\lambda}}
            \beta_{q_1} \cdots \beta_{q_{r-2\lambda}}
     \sum_{\delta_1|\pi_1(q)}
            \Bigl(
	      \frac{\pm n\delta_1 p}{\pi_2}
	    \Bigr)
       \frac{1}{\sqrt{\pi_1'}}
\Bigr]
                    \nonumber
\end{eqnarray}
Note that 
$$
n \delta_1 p \le n\pi_1 \le nU.
$$
By induction, each of the two inner $q$-sums is
$$
\ll
(2\epsilon_E e)^{r-2\lambda}
\Bigl(
  \log N_E + 2\log(|n|U) + \log (|n| x) + \log(U+2)
\Bigr)^{r-2\lambda}
\log^{2(r-2\lambda)+1}x.
$$
Also,
$$
\sum_{\lambda=1}^{\lfloor r/2 \rfloor}
\sum_{p<x}
      \Bigl(
        \frac{4\log^2 p}{p}
      \Bigr)^\lambda
\ll
\log^2 x
+
\sum_{\lambda\ge 2}
\sum_p 
      \Bigl(
        \frac{4\log^2 p}{p}
      \Bigr)^\lambda
\ll
\log^2 x.
$$
Thus (\ref{here}) is
$$
\ll
(2\epsilon_E e)^{r-2}
\Bigl(
  \log N_E + 2\log(|n|U) + \log (|n| x) + \log(U+2)
\Bigr)^{r-2}
\log^{2r-1}x.
$$
Combine this with (\ref{pione}) and recall that $U\le x^r$, we are done.

\if 3\
{
\vsp

\section{Large sieve inequality}
    \label{sec:sieve}

From the large sieve Gallanger \cite{gallagher} derives the following
inequality: given any sequence of complex numbers
$
\alpha_1, \alpha_2, \ldots
$,
\begin{equation}
  \sum_{q\le Q} \frac{q}{\phi(q)}
  \sideset{}{_{}^\ast}\sum_\chi
       \:
       \Bigl| \sum_{n\le x} \alpha_n \chi(n) \Bigr|^2
\le
  (x - 1 + Q^2)\sum_{n\le x} |\alpha_n|^2,          \label{sieve}
\end{equation}
where
$
\sum_\chi^\ast
$
denotes a sum over all primitive Dirichlet characters mod $q$.  Note that
\begin{eqnarray*}
\Bigl( \sum_{p<x} \beta_p \chi(p) \Bigr)^k
&=&
\sum_{n\le x^k}
  \chi(n)
  \sum_{\substack{(p_1, \cdots, p_k)\\p_1 \cdots p_k = n}}
    \beta_{p_1} \cdots \beta_{p_k},
\end{eqnarray*}
so
\begin{eqnarray}
&&
\sum_{q\le Q} \frac{q}{\phi(q)}
  \sideset{}{_{}^\ast}\sum_\chi
    \:
    \Bigl|
       \Bigl(
         \frac{2}{\log x}
         \sum_{p\le x}
           \beta_p \chi(p)
       \Bigr)^k
    \Bigr|^{2}
                     \nonumber
\\
&\le&
\frac{2^{2k}(x^k - 1 + Q^2)}{\log^{2k} x}
\sum_{n\le x^k}
  \:\:\:
  \Bigl|
  \!\!\!
  \sum_{\substack{(p_1, \ldots, p_k)\\p_1 \cdots p_k = n}}
    \beta_{p_1} \cdots \beta_{p_k}
  \Bigr|^2
                     \nonumber
\\
&\le&
\frac{2^{2k}(x^k - 1 + Q^2)}{\log^{2k} x}
\sum_{n\le x^k}
  \sum_{\substack{(p_1, \ldots, p_k)\\p_1 \cdots p_k = n}}
    \beta_{p_1}^2 \cdots \beta_{p_k}^2
\hspace{20pt}
\text{unique factorization}
                     \nonumber
\\
&\le&
\frac{2^{2k}(x^k - 1 + Q^2)}{\log^{2k} x}
\Bigl(
  \sum_{p\le x}  \beta_p^2
\Bigr)^k
\hspace{84pt}
\text{since $\beta_p=0$ if $p>x$.}
                     \label{since}
\end{eqnarray}
Define
$$
P(k,x)
=
\{ \mymod{\chi}{q}:
   \text{$\chi$ primitive, $q \le x^{k/2}\log x$} \}.
$$
For a fixed $\epsilon>0$, and suppose there exists a constant
$
C_k(\epsilon)>0
$
so that
\begin{eqnarray}
\lefteqn{\#\Bigl\{
  \chi\in P(k,x):
      \text{$ \Bigl|\frac{2}{\log x}\sum_{p\le x} \beta_p \chi(p)\Bigr|
	   >
	  1/\sqrt{3} + \epsilon$}
\Bigr\}}
          \hspace{50pt}
          \nonumber
\\
&\ge&
C_k(\epsilon) \#P(k, x)
\:\:=\:\:
C_k(\epsilon) \Bigl( \frac{3}{\pi^2} + o(1) \Bigr) x^k\log^2 x
              \label{set}
\end{eqnarray}
(we are \textit{not} assuming that the set in question has a density).
Combine this with (\ref{since}) and we get
$$
C_k(\epsilon)
\le
\frac{\pi^2}{3}
\Bigl(
  \frac{1/\sqrt{3}}{1/\sqrt{3}+\epsilon}
\Bigr)^k
=
\frac{\pi^2}{3}
\Bigl(
  \frac{1}{1+\sqrt{3}\epsilon}
\Bigr)^k.
$$
The number of $\chi\in P(k,x)$ of conductor $\le x^{k/2}$ is
$
O(x^k)
$,
and if
$
q \ge x^{k/2}
$
then
$
2\log q = (k+o(1))\log x
$.
Recall the twisted explicit formula (\ref{twisted}) and we see that, for
$k\ge 1$, the set of
$
\chi\in P(k, x)
$
such that
$$
2k + 1/2 - 1/\sqrt{3} - \epsilon
\le
           \sum_{\tau_\chi} \Phi_{\log x}(\rho_\chi)
\le
2k + 1/2 + 1/\sqrt{3} + \epsilon
$$
has density
$$
\le  \frac{\pi^2}{3}
\Bigl(
  \frac{1}{1+\sqrt{3}\epsilon}
\Bigr)^k,
$$
as desired.
}
\fi

\vtsp

\bibliographystyle{amsalpha}

\end{document}